\documentclass[12pt]{article}%
\usepackage{amsmath}
\usepackage{amsfonts}
\usepackage{amssymb}
\usepackage{graphicx}
\usepackage{graphics}%
\newtheorem{theorem}{Theorem}
\newtheorem{acknowledgement}[theorem]{Acknowledgement}

\newtheorem{result}[theorem]{Result}

\begin{document}

\title{Geometrical Optics Approach to Markov-Modulated Fluid Models}
\author{Diego Dominici\thanks{ddomin1@uic.edu} \ and Charles
Knessl\thanks{knessl@uic.edu}\\Department of Mathematics, Statistics and Computer Science \\University of Illinois at Chicago (M/C 249) \\851 South Morgan Street \\Chicago, IL 60607-7045, USA}
\maketitle

\begin{abstract}
We analyze asymptotically a differential-difference equation, that arises in a
Markov-modulated fluid model. Here there are $N$ identical sources that turn
\textbf{on} and \textbf{off}, and when \textbf{on} they generate fluid at unit
rate into a buffer, which process the fluid at a rate $c<N.$ In the steady
state limit, the joint probability distribution of the buffer content and the
number of active sources satisfies a system of $N+1$ ODEs, that can also be
viewed as a differential-difference equation analogous to a backward/forward
parabolic PDE. We use singular perturbation methods to analyze the problem for
$N\rightarrow\infty,$ with appropriate scalings of the two state variables. In
particular, the ray method and asymptotic matching are used.

\end{abstract}

\section{Introduction}

The study of fluid queues has been the subject of much recent work. In these
models the queue length is considered a deterministic (or \textquotedblleft
fluid\textquotedblright) process, rather than a discrete random process that
measures the number of customers. These models tend to be somewhat easier to
analyze, as they allow for less randomness than more traditional queueing
models. Also, the rougher description of a queue as a fluid is thought to be
adequate for many important modern applications, such as ATM (asynchronous
transfer mode) and other high speed integrated networks.

Some of the earliest studies of fluid queues are due to Kosten \cite{Kosten1},
\cite{Kosten2} and Anick, Mitra and Sondhi \cite{AMS}. We briefly describe the
model in \cite{AMS}, since much of the latter work can be viewed as
extensions/generalizations of it. There are $N$ identical, independent sources
that turn \textbf{on} and \textbf{off} at exponential waiting times. The rate
from \textbf{off} to \textbf{on} is $\lambda,$ and time is scaled so that the
rate of the reverse transition equals $1.$ When a particular source is
\textbf{on}, it generates fluid at unit rate to a buffer. We denote the buffer
content at time $t$ by $X(t).$ This buffer is depleted at a constant
(deterministic) rate $c,$ as long as it is non-empty.

If $c<N,$ which we assume henceforth, the buffer content may be non-empty.
Letting $Z(t)$ be the number of \textbf{on} sources at time $t,$ the buffer
evolves according to the law
\[
\dot{X}(t)=\left\{
\begin{array}
[c]{c}%
Z(t)-c,\text{ if}\quad X(t)>0\\
\max\left[  Z(t)-c,\ 0\right]  ,\text{ if}\quad X(t)=0
\end{array}
\right.
\]
where $\left\lfloor c\right\rfloor $ denotes the integer part of $c$. For
simplicity we assume that $c$ is not an integer, so that the fractional part
of $c,$ $\left\{  c\right\}  =c-\left\lfloor c\right\rfloor \in(0,1).$

Since there is randomness in the evolution of $Z(t),$ the process
$(X(t),Z(t))$ is stochastic and the exponential distributions of the
\textbf{on} and \textbf{off} periods imply that the process is Markovian.
These models are called Markov-modulated fluid flows. Of primary interest in
applications is the distribution of the buffer content $X(t).$ This can be
obtained easily if we know the joint distribution $F_{k}(x,t)$%
\begin{align*}
F_{k}(x,t)  &  =\Pr\left[  X(t)\leq x,\ Z(t)=k\mid X(0)=x_{0},\ Z(0)=k_{0}%
\right] \\
0  &  \leq x,x_{0}\ ;\ 0\leq k,k_{0}\leq N
\end{align*}
but the latter is often difficult to compute. It seems that it is known
explicitly only for the simplest models, such as those in \cite{Kosten1},
\cite{Kosten2} and \cite{AMS}. Surveys of (single buffer) fluid models appear
in \cite{Kultrami}, \cite{Resing} and \cite{Mitra}.

In \cite{AMS} the steady state limit $t\rightarrow\infty$ is considered for
the model described above and the authors obtained an explicit expression for
$F_{k}(x)=F_{k}(x,\infty),$ as a spectral representation. Then, the marginal
buffer distribution was obtained from $%
{\displaystyle\sum\nolimits_{k=0}^{N}}
F_{k}(x).$ The complexity of this function has led to various asymptotic
investigations. These assume that the number of sources $N\rightarrow\infty,$
with an appropriate scaling of the state variables $(x,k).$ In \cite{S Weiss}
and \cite{Weiss1} the theory of large deviations is used to construct an
approximation to the probability that the buffer exceeds the value $x=Ny,$ of
the form $\Pr\left[  X(\infty)>Ny\right]  \approx\exp\left[  -N\ \mathrm{I}%
(y)\right]  ,$ where $\mathrm{I}(y)$ satisfies a variational problem that can
be solved fairly explicitly.

A more complete asymptotic description of this buffer overflow probability is
obtained in \cite{Morrison}, where the author shows that
\[
\Pr\left[  X(\infty)>Ny\right]  \sim N^{-\frac{1}{2}}\mathrm{I}_{1}%
(y)\exp\left[  -N\ \mathrm{I}_{0}(y)\right]
\]
and it is also shown that $\mathrm{I}_{0}(y)$ is equivalent to the solution of
the variational problem in \cite{Weiss1}. It is furthermore established that
the asymptotic approximation is quite accurate numerically. The technique used
is to expand the spectral representation in \cite{AMS} using asymptotic
methods, such as the Euler-McLaurin formula and Laplace's method
\cite{BenderOrzag}.

The existence of a steady state for this model requires that $\frac{\lambda
}{\lambda+1}<\gamma<1,\quad c=N\gamma,$ which says simply that the processing
rate exceeds the average input rate. The case where this stability condition
is only weakly satisfied is called \textquotedblleft heavy
traffic\textquotedblright. More precisely, this corresponds to the scaling
$\gamma=\frac{\lambda}{\lambda+1}+O\left(  N^{-1/2}\right)  .$ In
\cite{Knessl-Morrison} the authors obtain an approximate diffusion model in
the heavy traffic limit and solve it explicitly as a spectral representation
involving Hermite polynomials. We note that the discrete model has
$F_{k}(0)\neq0$ for $0\leq k\leq\left\lfloor c\right\rfloor ,$ which says that
there is a non-zero probability that the buffer is empty for this range of
\textbf{on} sources. Then $F_{k}(0)=0$ for $k>c$ and this boundary condition
was used to obtain explicitly the coefficients in the spectral representation
in \cite{AMS}.

In the heavy traffic limit the joint distribution satisfies a parabolic PDE,
that behaves as the heat equation in a part of the domain and the backward
heat equation in the remainder \cite{Knessl-Keller}, \cite{Knessl-Morrison}.
Such problems arise in a variety of applications, such as counter-current
separators \cite{Hagan-Okendon}, mean exit times \cite{HDL}, the Milne problem
of statistical physics \cite{Bardos}, neutron transport theory \cite{Hopf} and
diffusion in spatially varying convection fields.

Their study goes back to Gevrey \cite{Gevrey1}, \cite{Gevrey 2} and more
recent analyses appear in \cite{Aziz Liu}, \cite{Aziz 2}, \cite{Baquendi},
\cite{Freidlin weinb} and \cite{Keller Wainb}. The interesting mathematical
feature of these problems is that the initial (or boundary) conditions can be
imposed only where the PDE is forward parabolic. This \textquotedblleft
half-boundary condition\textquotedblright\ makes the problem difficult to
analyze. The model in \cite{AMS} may be viewed as a discrete analog of these
PDEs and the \textquotedblleft half BC\textquotedblright\ corresponds to the
condition $F_{k}(0)=0,$ that can be applied only for $\left\lfloor
c\right\rfloor +1\leq k\leq N.$

In \cite{Knessl-Keller} we developed an asymptotic approach to analyze
backward-forward parabolic PDEs, in a limit where the diffusion coefficient is
small. It is based on the ray method of geometrical optics \cite{Keller} and
matched asymptotic expansions. An important feature is the careful treatment
of the point on the boundary where the PDE changes type from backward to
forward parabolic. In this \textquotedblleft corner region\textquotedblright%
\ the asymptotic solution may be represented as a contour integral involving
Airy functions.

The purpose of this paper is to extend the asymptotic approach in
\cite{Knessl-Keller} to discrete models. We shall analyze the model in
\cite{AMS} directly by using the differential-difference equation satisfied by
$F_{k}(x).$ We make no recourse to the exact spectral representation of the
solution given in \cite{AMS}. After appropriate scalings of $k$ and $x,$ we
shall analyze this equation asymptotically for $N\rightarrow\infty$ using
singular perturbation methods. The asymptotic results for $%
{\displaystyle\sum\nolimits_{k=0}^{N}}
F_{k}(x)$ in \cite{Morrison} can be easily recovered from our two-dimensional
results. There are several important differences between the asymptotic
structure of these discrete models and the backward-forward parabolic PDE
studied in \cite{Knessl-Keller}. For example, the structure of the solution in
the corner region, where $x\approx0$ and $k\approx c,$ is different than the
corresponding range in the diffusion model. Here the solution can be expressed
in terms of Bessel functions.

A variety of extensions/generalizations of the model in \cite{AMS} have been
recently analyzed. For example, the transient solution (more precisely, its
Laplace transform over time) is obtained in \cite{Kobashashi} as a spectral
expansion, using arguments similar to those in \cite{AMS}. In \cite{Maier} the
author allows the fluid input rates from the \textbf{on} sources to depend
upon the buffer size. This allows for an admission control policy for the
fluid level (buffer size). An asymptotic analysis was done in \cite{Maier},
which assumed that $N$ is fixed, but that the input rates vary
\textquotedblleft weakly\textquotedblright\ with the buffer size $x$. Finite
buffer size models are considered in \cite{Mitra} and in \cite{Mitra 18} the
model allows for certain sources to increase the buffer while others deplete
it. The latter leads to a three dimensional problem, where one must keep track
of the buffer content and also the number of active sources of either type.
Problems with two buffers and various priority mechanisms are studied in
\cite{Choi Choi 20}, \cite{Elwalis Mitra 19} and \cite{Liu Gong 21}. In
\cite{Doorn Sche 22} and \cite{Doorn Sche 23} fluid models with more general
birth-death modulating processes were analyzed.

It seems that for the more general models considered \cite{Choi Choi 20},
\cite{Doorn Sche 22}, \cite{Doorn Sche 23}, \cite{Elwalis Mitra 19}, \cite{Liu
Gong 21}, \cite{Mitra 18}, the solutions are not particularly explicit. They
involve either solving systems of equations or computing the eigenvalues of
matrices numerically. The merit of our asymptotic approach is that it yields
relatively simple formulas, which are both easy to numerically evaluate and
also provide numerically accurate approximations to the performance measures,
even for moderate values of the large parameter $N.$ We use similar scalings
as the large deviations studies \cite{Mandjes Ridder 24}, \cite{S Weiss},
\cite{Weiss1}, but in contrast to these studies we provide the full asymptotic
approximation and not just the exponential growth/decay rate (in $N)$ of the
performance measure.

We also carefully treat various boundary and corner regions of the state
space
\begin{equation}
\left\{  \left(  x,k\right)  :x\geq0,\quad0\leq k\leq N\right\}  \label{SS}%
\end{equation}
and indeed we show that their analysis is needed in order to obtain the
asymptotic expansions away from the boundaries. We obtain detailed results for
the model in \cite{AMS} and develop the methodology to treat other models of
this type.

The paper is organized as follows. In Section 2 we state the basic equations.
In Sections 3-7 we analyze these in various ranges of the state space
(\ref{SS}). In Section 8\ we recover the one-dimensional results in
\cite{Morrison}. Finally, in Section 9 we summarize and interpret the asymptotics.

\section{Problem statement}

In the model proposed by Amick, Mitra and Sondhi \cite{AMS}, a data-handling
switch receives messages from $N$ mutually independent information sources,
which independently and asynchronously alternate between the \textbf{on}\ and
\textbf{off}\ state. The number of \textbf{on}\ sources forms a birth-death
process $Z(t)$ with birth rate $\lambda_{k}=\lambda(N-k)$ and death rate
$\mu_{k}=k,$ where the rates are conditioned on $Z(t)=k$. Each source is
\textbf{on}, on average, $\frac{\lambda}{\lambda+1}$ of the total time. An
\textbf{on}\ source transmits at the uniform rate of $1$ unit of information
per unit of time.

The switch has infinite capacity, and stores or buffers the incoming
information that is in excess of the maximum transmission rate $c$\ of the
output channel. The drift $r_{k}=k-c$ gives the rate of increase of $X(t)$
(the buffer content at time $t)$ when the birth-death process is in state $k.$
That is, the rate of change of $X(t)$ at time $t$ is $r_{Z(t)},$ provided
$r_{Z(t)}\geq0$ or $r_{Z(t)}<0$ and $X(t)>0.$ If the buffer has emptied at
time $t$, it stays empty as long as the drift remains negative.

Following \cite{Scheinhardt} we define
\[
\pi_{k}=%
{\displaystyle\prod\limits_{j=0}^{k-1}}
\frac{\lambda_{j}}{\mu_{j+1}}=\binom{N}{k}\lambda^{k}.
\]
The stationary probabilities $p_{k}$ of the birth-death process can then be
represented as%
\[
p_{k}=\frac{\pi_{k}}{%
{\displaystyle\sum\nolimits_{j=0}^{N}}
\pi_{j}}=\frac{1}{\left(  \lambda+1\right)  ^{N}}\binom{N}{k}\lambda^{k}.
\]
In order that a stationary distribution for $X(t)$ exists, the mean drift $%
{\displaystyle\sum\nolimits_{j=0}^{N}}
\pi_{j}r_{j}$ should be negative%
\[
\left(  \lambda+1\right)  ^{N}\left[  \frac{\lambda}{\lambda+1}N-c\right]  =%
{\displaystyle\sum\nolimits_{j=0}^{N}}
\pi_{j}r_{j}<0
\]
which gives the stability condition%
\begin{equation}
\frac{\lambda}{\lambda+1}<\gamma<1,\quad\gamma=\frac{c}{N}. \label{stability}%
\end{equation}
Setting
\[
F_{k}(x,t)=\Pr\left[  X(t)\leq x,\ Z(t)=k\right]  ;\quad t,x\geq0,\quad0\leq
k\leq N
\]
and
\[
F_{k}(x,t)\equiv0,\quad k\notin\lbrack0,\ N],
\]
the Kolmogorov forward equations for the Markov process $\left(
X(t),\ Z(t)\right)  $ are given by%
\[
\frac{\partial F_{k}}{\partial t}+r_{k}\frac{\partial F_{k}}{\partial
x}=\lambda_{k-1}F_{k-1}+\mu_{k+1}F_{k+1}-\left(  \lambda_{k}+\mu_{k}\right)
F_{k}.
\]
For the stationary distribution $F_{k}(x)\equiv$ $F_{k}(x,\infty)$ with the
above rates and drift, we have%
\begin{equation}
(k-c)\frac{\partial F_{k}}{\partial x}=\lambda\left[  N-(k-1)\right]
\,\,F_{k-1}+(k+1)\,\,F_{k+1}-\left[  \lambda(N-k)+k\right]  \,F_{k},\quad0\leq
k\leq N. \label{diffeq}%
\end{equation}
Moreover, if the number of \textbf{on} sources $k$ exceeds $c$, then the
buffer content increases and the buffer can't be empty. Hence,%
\begin{equation}
F_{k}(0)=0,\quad\left\lfloor c\right\rfloor +1\leq k\leq N. \label{FBC}%
\end{equation}
Also,%
\begin{equation}
F_{k}(\infty)=\frac{1}{(1+\lambda)^{N}}\binom{N}{k}\lambda^{k},\quad0\leq
k\leq N, \label{Finf}%
\end{equation}
since $F_{k}(\infty)$ is the probability that $k$ out of $N$ sources are
\textbf{on}\ simultaneously.

\section{The ray expansion}

To analyze the problem (\ref{diffeq})-(\ref{Finf}) for large \ $N$ \ we
introduce the scaled variables $y,z,\gamma,$ with
\[
k=zN,\quad c=\gamma N,\quad x=yN,\quad z,\gamma,y=O(1).
\]
We define the function \ $G(y,z)$ \ and the small parameter \ $\varepsilon$
\ by \ \
\[
\varepsilon=\frac{1}{N},\quad F_{k}(x)=G\left(  \varepsilon x,\varepsilon
k\right)  =G(y,z)
\]
and note that \ $F_{k\pm1}(x)=G(y,z\pm\varepsilon).$

Then (\ref{diffeq}) becomes the following equation for\ $G(y,z)$\
\begin{equation}
\varepsilon(z-\gamma)G_{y}(y,z)=\lambda(1-z-\varepsilon)G(y,z-\varepsilon
)+(z+\varepsilon)G(y,z+\varepsilon)-\left[  \lambda(1-z)+z\right]  G(y,z)
\label{eqG}%
\end{equation}
and (\ref{FBC}) implies that%
\begin{equation}
G(0,z)=0,\quad\gamma<z<1. \label{BC}%
\end{equation}

To find \ $G(y,z)$ for\ \ $\varepsilon$ \ small, we shall use the ray method.
Thus, we consider solutions which have the asymptotic form
\begin{equation}
G(y,z)\sim\varepsilon^{\nu}\exp\left[  \frac{1}{\varepsilon}\Psi(y,z)\right]
\mathbb{K}(y,z). \label{GRay}%
\end{equation}
Using\ (\ref{GRay}) in (\ref{eqG}), with
\[
\frac{1}{\varepsilon}\Psi(y,z\pm\varepsilon)=\frac{1}{\varepsilon}\Psi\pm
\Psi_{z}+\frac{1}{2}\Psi_{zz}\varepsilon+O\left(  \varepsilon^{2}\right)  ,
\]
dividing by \ $\exp\left[  \frac{1}{\varepsilon}\Psi(y,z)\right]  ,$ and
expanding in powers of \ $\varepsilon$ \ we get
\begin{gather*}
(z-\gamma)\left[  \Psi_{y}\mathbb{K}+\varepsilon\mathbb{K}_{y}\right]
=\left[  zU+(\lambda-1)z-\lambda+\lambda(1-z)\frac{1}{U}\right]  \mathbb{K}\\
+\left\{  \left[  zU+\lambda(z-1)\frac{1}{U}\right]  \mathbb{K}_{z}+\left[
\left(  1+\frac{z}{2}\Psi_{zz}\right)  U+\lambda\left(  1+\frac{1}{2}\Psi
_{zz}-\frac{z}{2}\Psi_{zz}\right)  \frac{1}{U}\right]  \mathbb{K}\right\}
\varepsilon+O\left(  \varepsilon^{2}\right)
\end{gather*}
where%
\begin{equation}
U(y,z)=\exp\left[  \Psi_{z}(y,z)\right]  . \label{U}%
\end{equation}
Equating the coefficients of \ $\varepsilon$ \ we obtain the \emph{eikonal
equation} for \ $\Psi(y,z)$%
\begin{equation}
(z-\gamma)\Psi_{y}+(1-\lambda)z+\lambda+\lambda(z-1)\frac{1}{U}-zU=0
\label{eik}%
\end{equation}
and the \emph{transport equation} for $\mathbb{K}(y,z)$%
\begin{equation}
\left[  \left(  1+\frac{z}{2}\Psi_{zz}\right)  U+\lambda\left(  1+\frac{1}%
{2}\Psi_{zz}-\frac{z}{2}\Psi_{zz}\right)  \frac{1}{U}\right]  \mathbb{K}%
+(\gamma-z)\mathbb{K}_{y}+\left[  zU+\lambda(z-1)\frac{1}{U}\right]
\mathbb{K}_{z}=0. \label{trans}%
\end{equation}
To solve (\ref{eik}) and (\ref{trans}) we use the method of characteristics,
which we briefly review below.

Given the first order partial differential equation%
\[
\mathfrak{F}\left(  y,z,\Psi,p,q\right)  =0,
\]
where \ $p=\Psi_{y},\quad q=\Psi_{z},$ we search for a solution \ $\Psi(y,z).$
The technique is to solve the system of \textquotedblleft characteristic
equations\textquotedblright\ given by%
\begin{align*}
\dot{y}  &  =\frac{\partial y}{\partial t}=\mathfrak{F}_{p},\quad\dot
{z}=\mathfrak{F}_{q}\\
\dot{p}  &  =-\mathfrak{F}_{y}-p\mathfrak{F}_{\Psi},\quad\dot{q}%
=-\mathfrak{F}_{z}-q\mathfrak{F}_{\Psi}\\
\dot{\psi}  &  =p\mathfrak{F}_{p}+q\mathfrak{F}_{q}%
\end{align*}
where we now consider $\left\{  y,z,\psi,p,q\right\}  $ to all be functions of
the variable $t,$ with $\psi(s,t)=\Psi(y,z).$

For the eikonal equation (\ref{eik}), the characteristic equations are
\begin{subequations}
\label{strip}%
\begin{align}
\dot{y}  &  =z-\gamma\label{eqa}\\
\quad\dot{z}  &  =\lambda(1-z)e^{-q}-ze^{q}\label{eqb}\\
\dot{p}  &  =0\label{eqc}\\
\quad\dot{q}  &  =e^{q}-\lambda e^{-q}-p+\lambda-1\label{eqd}\\
\dot{\psi}  &  =p(z-\gamma)+q\left[  \lambda(1-z)e^{-q}-ze^{q}\right]  .
\label{eqe}%
\end{align}
The particular solution is determined by the initial conditions at $t=0$. We
shall show that for this problem two different types of solutions are needed;
these correspond to two distinct families of characteristic curves, or rays.

\subsection{The partial derivatives $\Psi_{y},\ \Psi_{z}$}

Setting \ $\left.  \Psi_{y}\right\vert _{t=0}=s$, \ and solving (\ref{eqc})
yields%
\end{subequations}
\[
p=s,
\]
so that $\Psi_{y}$ is constant along a ray. Introducing the function
\ $u(s,t)$ \ as in (\ref{U})%
\[
u(s,t)=\exp[q(s,t)]
\]
we have from (\ref{eqd})%
\begin{align}
\dot{u}  &  =u\dot{q}=\left(  u-\frac{\lambda}{u}-s+\lambda-1\right)
u\nonumber\\
&  =u^{2}+\left(  \lambda-1-s\right)  u-\lambda\label{u1}\\
&  =(u-r_{1})(u-r_{2})\nonumber
\end{align}
where%
\[
r_{1,2}(s)=\frac{1}{2}\left(  s+1-\lambda\pm\Delta\right)  ,\quad
\Delta(s)=\sqrt{\left(  \lambda-s-1\right)  ^{2}+4\lambda}%
\]
and $r_{1}$ corresponds to the $(+)$ sign.

Solving (\ref{u1}) gives%
\begin{equation}
\frac{1}{\Delta}\ln\left[  \frac{\left(  u-r_{1}\right)  \left(  u_{0}%
-r_{2}\right)  }{\left(  u-r_{2}\right)  \left(  u_{0}-r_{1}\right)  }\right]
=t,\quad u(s,0)=u_{0}(s) \label{t(u)}%
\end{equation}
or%
\begin{equation}
u(s,t)=r_{2}+\frac{\Delta}{e^{\Delta t}\left(  \frac{\Delta}{u_{0}-r_{2}%
}-1\right)  +1}. \label{u(s,t)}%
\end{equation}
Evaluating (\ref{eik}) at $t=0$ we get%
\[
(1-\lambda)\gamma+\lambda+\lambda(\gamma-1)\frac{1}{u_{0}}-\gamma u_{0}=0
\]
so that%
\begin{equation}
u_{0}=1\text{ \ or \ }u_{0}=\frac{\lambda}{\gamma}\left(  1-\gamma\right)  .
\label{u0}%
\end{equation}

\subsection{The rays from $(0,\gamma)$}

We now consider the family of rays emanating from the point \ $y=0,\ z=\gamma
.$ \ From (\ref{eqb}) the equation for \ $z$ \ is%
\begin{equation}
\dot{z}=\lambda(1-z)\frac{1}{u}-zu \label{dzdt}%
\end{equation}
which we can rewrite as%
\[
\frac{dz}{du}\dot{u}=\lambda(1-z)\frac{1}{u}-zu
\]
whose solution is$\allowbreak$%
\[
z=\frac{C(s)u-\lambda}{u^{2}+\left(  \lambda-1-s\right)  u-\lambda}.
\]
Imposing the initial condition \ $z(s,0)=\gamma$ and using (\ref{u0}), \ we
obtain $C(s)=\lambda-\gamma s$ for both possible values of $u_{0}.$ Therefore,%
\begin{equation}
z=\frac{(\lambda-\gamma s)u-\lambda}{u^{2}+\left(  \lambda-1-s\right)
u-\lambda}. \label{z(u)}%
\end{equation}

From (\ref{strip}) the equation for \ $y$ \ is%
\begin{equation}
\dot{y}=z-\gamma\label{dydt}%
\end{equation}
which implies that \
\[
\dot{y}(s,0)=z(s,0)-\gamma=0
\]
and%
\[
\ddot{y}(s,0)=\dot{z}(s,0)=\lambda(1-\gamma)\frac{1}{u_{0}}-\gamma u_{0}.
\]
We define
\begin{equation}
\rho=\gamma-\lambda+\lambda\gamma\label{rho}%
\end{equation}
with $0<\rho<1$ from (\ref{stability}). From (\ref{u0}) we have%
\[
\ddot{y}(s,0)=\left\{
\begin{array}
[c]{l}%
-\rho,\quad u_{0}=1\\
\allowbreak\rho,\quad u_{0}=\frac{\lambda}{\gamma}\left(  1-\gamma\right)
\end{array}
\right.  .
\]
Using the initial condition \ $y(s,0)=0$ \ and expanding in powers of \ $t$,
\ we get%
\[
y(s,t)\sim\ddot{y}(s,0)\frac{t^{2}}{2},\quad t\rightarrow0
\]
and in order to have \ $y>0$ \ for $t>0$ (i.e., for the rays to enter the
domain $[0,\infty)\times\lbrack0,1]$) we need to choose
\begin{equation}
\ u_{0}=\frac{\lambda}{\gamma}\left(  1-\gamma\right)  \label{U0}%
\end{equation}
with $\ u_{0}<1$ from (\ref{stability}).

Integrating (\ref{dydt}) and using (\ref{u1}), (\ref{t(u)}) and (\ref{z(u)}),
we have%
\begin{align*}
y  &  =%
{\displaystyle\int\limits_{0}^{t}}
z(s,v)dv-\gamma t\\
&  =%
{\displaystyle\int\limits_{u_{0}}^{u}}
\frac{(\lambda-\gamma s)w-\lambda}{\left[  w^{2}+\left(  \lambda-1-s\right)
w-\lambda\right]  ^{2}}dw-\gamma\frac{1}{\Delta}\ln\left[  \frac{\left(
u-r_{1}\right)  \left(  u_{0}-r_{2}\right)  }{\left(  u-r_{2}\right)  \left(
u_{0}-r_{1}\right)  }\right]
\end{align*}%
\begin{align}
&  =-\frac{1}{\Delta^{3}}\left[  (\lambda+1)\rho+\phi s\right]  \ln\left[
\frac{\left(  u-r_{1}\right)  \left(  u_{0}-r_{2}\right)  }{\left(
u-r_{2}\right)  \left(  u_{0}-r_{1}\right)  }\right] \label{y(u)}\\
&  +\frac{\left[  r_{1}(\lambda-s\gamma)-\lambda\right]  (u-u_{0})}{\Delta
^{2}\left(  u-r_{1}\right)  \left(  u_{0}-r_{1}\right)  }+\frac{\left[
r_{2}(\lambda-s\gamma)-\lambda\right]  (u-u_{0})}{\Delta^{2}\left(
u-r_{2}\right)  \left(  u_{0}-r_{2}\right)  }\nonumber
\end{align}
where%
\begin{equation}
\phi=\gamma+\lambda-\gamma\lambda\label{phi}%
\end{equation}
with $\phi>0$ from (\ref{stability}).

Finally, combining (\ref{u(s,t)}), (\ref{z(u)}), (\ref{y(u)}) and (\ref{U0})
we conclude that%
\begin{equation}
y(s,t)=\frac{1}{\Delta^{2}}\left\{  \frac{\phi s+(\lambda+1)\rho}{\Delta}%
\sinh(\Delta t)+\rho\cosh(\Delta t)-\rho+\left[  \gamma s^{2}+(\gamma
-\lambda-\lambda\gamma)s+\lambda(\lambda+1)\right]  t\right\}  -\gamma t
\label{y(s,t)}%
\end{equation}%
\begin{equation}
z(s,t)=\frac{1}{\Delta^{2}}\left[  \phi s+(\lambda+1)\rho\cosh(\Delta
t)+\rho\Delta\sinh(\Delta t)+\gamma s^{2}+(\gamma-\lambda-\lambda
\gamma)s+\lambda(\lambda+1)\right]  . \label{z(s,t)}%
\end{equation}
This yields the rays that emanate from $(0,\gamma)$ in parametric form.
Several rays are sketched in Figure 1.

\begin{figure}[t]
\begin{center}
\rotatebox{270} {\resizebox{10cm}{!}{\includegraphics{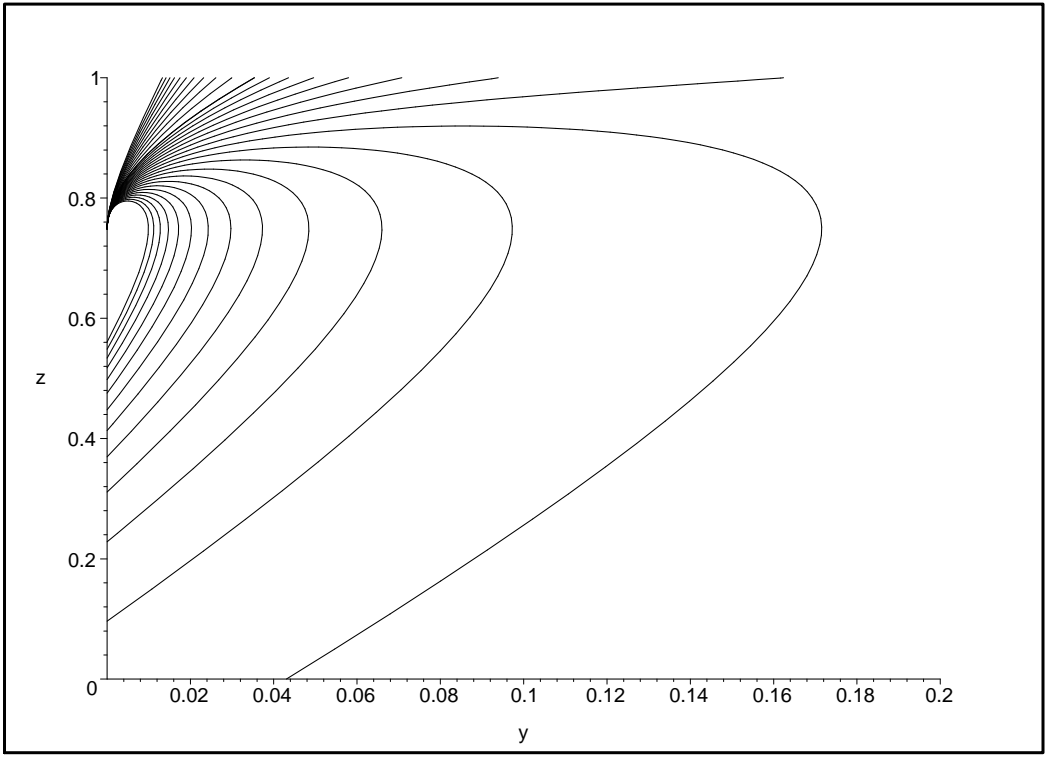}}}
\end{center}
\caption{A sketch of the rays from $(0,\gamma)$.}%
\end{figure}

Combining the terms depending on \ $t$ \ we have another useful expression for
\ $y(s,t)$ \ and \ $z(s,t)$%
\begin{equation}
y(s,t)=\frac{A_{1}(s)}{\Delta}\left(  e^{\Delta t}-1\right)  +\frac{A_{2}%
(s)}{\Delta}\left(  1-e^{-\Delta t}\right)  +\left[  A_{3}(s)-\gamma\right]  t
\label{y(t)}%
\end{equation}%
\begin{equation}
z(s,t)=A_{1}(s)e^{\Delta t}+A_{2}(s)e^{-\Delta t}+A_{3}(s) \label{z(t)}%
\end{equation}
with%
\begin{align*}
A_{1}(s)  &  =\frac{1}{4\Delta^{2}}\left(  \lambda+1-s+\Delta\right)  \left(
s\gamma+\rho-\lambda+\gamma\Delta\right) \\
A_{2}(s)  &  =\frac{1}{4\Delta^{2}}\left(  \lambda+1-s-\Delta\right)  \left(
s\gamma+\rho-\lambda-\gamma\Delta\right) \\
A_{3}(s)  &  =\frac{1}{\Delta^{2}}\left[  \lambda(1+\lambda)+\left(
\gamma-\gamma\lambda-\lambda\right)  s+\gamma s^{2}\right]  .
\end{align*}
\bigskip

For \ $t\geq0$ \ and each value of \ $s$, \ (\ref{y(s,t)}) and (\ref{z(s,t)})
determine a ray in the\ $(y,z)$ plane, which starts from \ $(0,\gamma)$ \ at
\ $t=0$. \ We discuss two particular rays which can be obtained in an explicit
form. For \ $s=0$ \ we can eliminate $\ t$ \ from (\ref{z(s,t)}) and obtain%
\begin{equation}
t=\frac{1}{\lambda+1}\ln\left(  \frac{z\lambda+z-\lambda}{\rho}\right)
,\quad\gamma<z<1 \label{T0}%
\end{equation}
and using (\ref{y(s,t)}) we get%
\begin{equation}
y=Y_{0}(z)=\frac{z-\gamma}{\lambda+1}-\frac{\rho}{\left(  \lambda+1\right)
^{2}}\ln\left(  \frac{z\lambda+z-\lambda}{\rho}\right)  ,\quad s=0,\quad
\gamma\leq z\leq1. \label{Y0}%
\end{equation}
Also, from (\ref{z(u)}) we have%
\begin{equation}
u=\lambda\frac{1-z}{z},\quad\text{for \ }s=0,\quad\gamma\leq z\leq1.
\label{u(0)}%
\end{equation}

For \ $s>S_{0}\equiv-\frac{\rho}{\gamma(1-\gamma)}$, \ we have both \ $y(t)$
\ and \ $z(t)$ \ increasing for \ $t>0$. When \ $s=S_{0}$ \ the function $u$
is constant along the ray%
\[
u(S_{0},t)\equiv u_{0}=\frac{\lambda}{\gamma}\left(  1-\gamma\right)  ,
\]
$y(t)$ \ increases and \ $z(t)$ \ asymptotes $\frac{\gamma^{2}}{\delta}$ with
\[
\delta=(\gamma-1)^{2}\lambda+\gamma^{2},\quad\delta>0.
\]
\ Eliminating \ $t$ \ we obtain \
\[
t=\frac{\gamma(1-\gamma)}{\delta}\ln\left[  \frac{\gamma(1-\gamma)\rho}%
{\gamma^{2}-\delta z}\right]  ,\quad\gamma\leq z<\frac{\gamma^{2}}{\delta}%
\]%
\begin{equation}
y=Y_{1}(z)=\frac{\left[  \gamma(1-\gamma)\right]  ^{2}\rho}{\delta^{2}}\left[
\frac{\gamma^{2}-\delta z}{\gamma(1-\gamma)\rho}+\ln\left[  \frac
{\gamma(1-\gamma)\rho}{\gamma^{2}-\delta z}\right]  -1\right]  ,\quad
s=S_{0},\quad\gamma\leq z<\frac{\gamma^{2}}{\delta}. \label{Y1}%
\end{equation}

For \ $s<S_{0}$ \ the rays reach a maximum value in \ $z$ \ at \ $t=T_{\max}$, where%

\[
T_{\max}=\frac{1}{2\Delta}\ln\left(  \frac{A_{2}}{A_{1}}\right)
\]
and we have%
\begin{equation}
y(s,T_{\max})=-\frac{\rho}{\Delta^{2}}+\frac{A_{3}-\gamma}{2\Delta}\ln\left(
\frac{A_{2}}{A_{1}}\right)  \label{ytmax}%
\end{equation}%
\begin{equation}
z(s,T_{\max})=2\sqrt{A_{1}A_{2}}+A_{3}. \label{ztmax}%
\end{equation}
Inverting (\ref{ztmax}) we find that at $t=T_{\max}$%
\begin{equation}
s=\frac{z+\lambda(1-z)\pm2\sqrt{z\lambda(1-z)}}{\gamma-z}. \label{smax}%
\end{equation}
From (\ref{ztmax}) we see that \ $z\rightarrow\frac{\gamma^{2}}{\delta}$ \ as
\ $s\uparrow S_{0}$. Therefore we must choose the $(-)$ sign in (\ref{smax})
and conclude that%
\begin{equation}
S_{\max}(z)=\frac{z+\lambda(1-z)-2\sqrt{z\lambda(1-z)}}{\gamma-z},\quad
\gamma<z<\frac{\gamma^{2}}{\delta} \label{Smax}%
\end{equation}
which when used in (\ref{ytmax}) yields%
\begin{equation}
y=Y_{\max}(z)=-\frac{\rho}{\Delta^{2}\left(  S_{\max}\right)  }+\frac
{A_{3}\left(  S_{\max}\right)  -\gamma}{2\Delta\left(  S_{\max}\right)  }%
\ln\left[  \frac{A_{2}\left(  S_{\max}\right)  }{A_{1}\left(  S_{\max}\right)
}\right]  ,\quad\gamma<z<\frac{\gamma^{2}}{\delta}. \label{Ymax}%
\end{equation}

\begin{figure}[t]
\begin{center}
\rotatebox{270} {\resizebox{10cm}{!}{\includegraphics{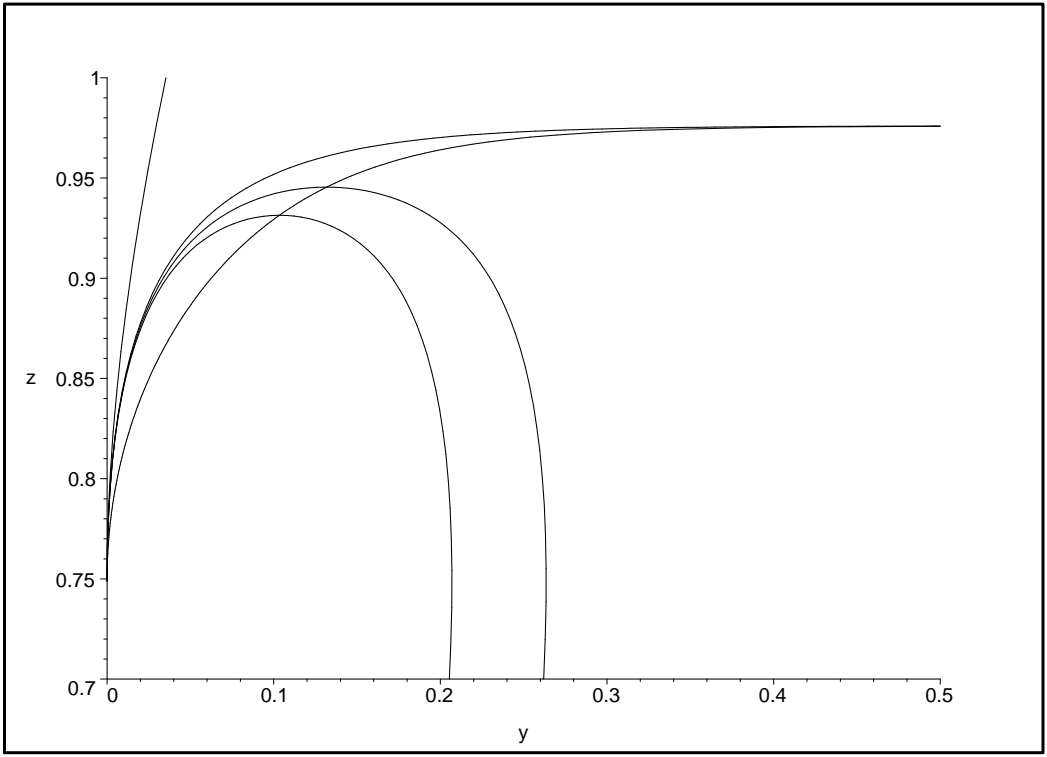}}}
\end{center}
\caption{A sketch of the rays $Y_{0}(z),Y_{1}(z)$ and the curve $Y_{\max}%
(z)$.}%
\end{figure}

This is the value of $y$ when $z$ reaches its maximum value. In Figure 2 we
sketch the rays $Y_{0}(z),$ $Y_{1}(z),$ the critical curve $Y_{\max}(z)$ and a
pair of rays showing their intersection with $Y_{\max}(z)$ at the maximum.

Solving for $t$ in (\ref{z(t)}) we obtain%
\begin{equation}
\mathbb{T(}s,z)=\left\{
\begin{array}
[c]{l}%
T_{+}(s,z),\quad0<y<Y_{\max}(z),\ \quad\gamma<z<1\\
T_{-}(s,z),\quad\left\{  y>Y_{\max}(z),\ \gamma<z<\frac{\gamma^{2}}{\delta
}\right\}  \cup\left\{  y>0,\ 0<z<\gamma\right\}
\end{array}
\right.  \label{T}%
\end{equation}
with%
\[
T_{\pm}(s,z)=\frac{1}{\Delta}\ln\left[  \frac{z-A_{3}\pm\sqrt{\left(
z-A_{3}\right)  ^{2}-4A_{1}A_{2}}}{2A_{1}}\right]  .
\]

Solving for \ $u$ \ in (\ref{z(u)}) we find that%
\begin{equation}
u=\left\{
\begin{array}
[c]{l}%
U_{-}(s,z),\quad0<y<Y_{\max}(z),\ \gamma<z<1\\
U_{+}(s,z),\quad\left\{  y>Y_{\max}(z),\ \gamma<z<\frac{\gamma^{2}}{\delta
}\right\}  \cup\left\{  y>0,\ 0<z<\gamma\right\}
\end{array}
\right.  \label{U+-}%
\end{equation}
where%
\begin{align*}
U_{\pm}(s,z)  &  =\frac{1}{2}\left[  s+1-\lambda+\frac{1}{z}\left(
\lambda-\gamma s\right)  \right] \\
&  \pm\frac{1}{2}\sqrt{\rho^{2}+\left[  2\left(  \lambda+1\right)  \rho+2\phi
s\right]  (z-\gamma)+\left[  (\lambda+1)^{2}+2(1-\lambda)s+s^{2}\right]
(z-\gamma)^{2}}.
\end{align*}

\bigskip\ From (\ref{dydt}) we see that the maximum value in $y$ is achieved
at the same time that \ $z=\gamma$, and that occurs at \ $t=T_{\gamma}$ \ with%
\begin{equation}
T_{\gamma}=\frac{1}{\Delta}\ln\left[  \frac{\phi s-\rho\Delta+\left(
\lambda+1\right)  \rho}{\phi s+\rho\Delta+\left(  \lambda+1\right)  \rho
}\right]  ,\quad s<S_{0}. \label{Tg}%
\end{equation}
When \ $z=\gamma,\ y>0$ \ we also have from (\ref{z(u)})%
\begin{equation}
u\left(  s,T_{\gamma}\right)  =1,\quad s<S_{0}, \label{ug}%
\end{equation}
so that $\Psi_{z}=0$ when $z=\gamma.$

Inverting the equations (\ref{y(s,t)})-(\ref{z(s,t)}) we can write%
\[
s=S\left(  y,z\right)  ,\quad t=T(y,z)
\]
and%
\[
\Psi(y,z)=\psi\left[  S\left(  y,z\right)  ,T(y,z)\right]  ,\quad
\mathbb{K(}y,z\mathbb{)}=K\left[  S\left(  y,z\right)  ,T(y,z)\right]  .
\]
We will use this notation in the rest of the article.

\subsection{The function \ $\Psi$}

From (\ref{eqe}) we have%
\[
\dot{\psi}=s\dot{y}+\ln(u)\dot{z}%
\]
which we can integrate to get%
\begin{align*}
\psi(s,t)  &  =\psi(s,0)+sy(s,t)+%
{\displaystyle\int\limits_{0}^{t}}
\ln\left[  u(s,r)\right]  \frac{dz}{dr}(s,r)dr\\
&  =\psi_{0}(s)+sy+\ln(u)z-\ln\left(  u_{0}\right)  \gamma-%
{\displaystyle\int\limits_{0}^{t}}
\frac{z}{u}du.
\end{align*}
Using (\ref{z(u)}) and (\ref{u1}) we can write%
\begin{align*}%
{\displaystyle\int\limits_{0}^{t}}
\frac{z}{u}du  &  =%
{\displaystyle\int\limits_{0}^{t}}
\frac{(\lambda-\gamma s)u-\lambda}{u(u-r_{1})(u-r_{2})}du=%
{\displaystyle\int\limits_{0}^{t}}
\left[  \frac{1}{u}+\frac{(\lambda-\gamma s)r_{1}-\lambda}{r_{1}%
(u-r_{1})\Delta}-\frac{(\lambda-\gamma s)r_{2}-\lambda}{r_{2}(u-r_{2})\Delta
}\right]  du\\
&  =\ln\left(  \frac{u}{u_{0}}\right)  +\frac{(\lambda-\gamma s)r_{1}-\lambda
}{r_{1}\Delta}\ln\left(  \frac{u-r_{1}}{u_{0}-r_{1}}\right)  -\frac
{(\lambda-\gamma s)r_{2}-\lambda}{r_{2}\Delta}\ln\left(  \frac{u-r_{2}}%
{u_{0}-r_{2}}\right)  .
\end{align*}
Hence,%

\begin{align}
\psi(s,t)  &  =\psi_{0}(s)+sy+(z-1)\ln(u)+(1-\gamma)\ln\left(  u_{0}\right)
+\frac{(\lambda-\gamma s)r_{2}-\lambda}{r_{2}\Delta}\ln\left(  \frac{u-r_{2}%
}{u_{0}-r_{2}}\right) \nonumber\\
&  -\frac{(\lambda-\gamma s)r_{1}-\lambda}{r_{1}\Delta}\ln\left(
\frac{u-r_{1}}{u_{0}-r_{1}}\right)  . \label{psi1}%
\end{align}
Obviously, $\psi(s,0)\equiv\psi_{0}$ is a constant since all rays start at the
same point. We will determine $\psi_{0}$ in section 3.

\subsection{The transport equation}

Re-writing the transport equation (\ref{trans}) as%
\[
\left[  \left(  1+\frac{z}{2}\Psi_{zz}\right)  U+\lambda\left(  1+\frac{1}%
{2}\Psi_{zz}-\frac{z}{2}\Psi_{zz}\right)  \frac{1}{U}\right]  \mathbb{K}%
=(z-\gamma)\mathbb{K}_{y}+\left[  zU+\lambda(z-1)\frac{1}{U}\right]
\mathbb{K}_{z}%
\]
and taking into consideration (\ref{dydt}) and (\ref{dzdt}), we see that%
\[
\left[  \left(  1+\frac{z}{2}\psi_{zz}\right)  u+\lambda\left(  1+\frac{1}%
{2}\psi_{zz}-\frac{z}{2}\psi_{zz}\right)  \frac{1}{u}\right]  K=\dot{y}%
K_{y}+\dot{z}K_{z}=\dot{K}%
\]
or%
\begin{align*}
\frac{\dot{K}}{K}  &  =\left(  1+\frac{1}{2}z\psi_{zz}\right)  u+\frac
{\lambda}{u}\left[  1+\frac{1}{2}(1-z)\psi_{zz}\right] \\
&  =u+\frac{\lambda}{u}+\frac{1}{2}\left[  zu+\frac{\lambda}{u}(1-z)\right]
\psi_{zz}\\
&  =\frac{1}{2}\left(  u+\frac{\lambda}{u}\right)  +\frac{1}{2}\frac{\partial
}{\partial z}\left[  zu-\frac{\lambda}{u}(1-z)\right]  ,
\end{align*}
since \ $u=\exp\left[  \frac{\partial\psi}{\partial z}\right]  .$ But,
\ $zu-\frac{\lambda}{u}(1-z)=-\dot{z},$ so that%
\begin{align*}
\frac{\dot{K}}{K}  &  =\frac{1}{2}\left(  u+\frac{\lambda}{u}\right)
-\frac{1}{2}\frac{\partial\dot{z}}{\partial z}\\
&  =\frac{1}{2}\left(  u+\frac{\lambda}{u}\right)  -\frac{1}{2}\left(
\ddot{z}\frac{\partial t}{\partial z}+\frac{\partial\dot{z}}{\partial s}%
\frac{\partial s}{\partial z}\right)  .
\end{align*}
Introducing the Jacobian of the transformation from Cartesian $(y,z)$ to ray
$(s,t)$ coordinates%
\begin{equation}
\mathbf{J}(s,t)=\frac{\partial z}{\partial t}\frac{\partial y}{\partial
s}-\frac{\partial z}{\partial s}\frac{\partial y}{\partial t} \label{J}%
\end{equation}
we have%
\[
\frac{\partial s}{\partial z}=-\frac{\dot{y}}{\mathbf{J}}=\frac{\gamma
-z}{\mathbf{J}},\quad\frac{\partial t}{\partial z}=\frac{1}{\mathbf{J}}%
\frac{\partial y}{\partial s}%
\]
and%
\[
\mathbf{\dot{J}}=\ddot{z}\frac{\partial y}{\partial s}+\dot{z}\frac
{\partial\dot{y}}{\partial s}-\frac{\partial\dot{z}}{\partial s}\dot{y}%
-\frac{\partial z}{\partial s}\ddot{y}=\ddot{z}\frac{\partial y}{\partial
s}-\frac{\partial\dot{z}}{\partial s}(z-\gamma).
\]
Hence,%
\begin{align*}
\frac{\dot{K}}{K}  &  =\frac{1}{2}\left(  u+\frac{\lambda}{u}\right)
-\frac{1}{2J}\left[  \ddot{z}\frac{\partial y}{\partial s}-\frac{\partial
\dot{z}}{\partial s}(z-\gamma)\right] \\
&  =\frac{1}{2}\left(  u+\frac{\lambda}{u}\right)  -\frac{\mathbf{\dot{J}}%
}{2J}=\frac{1}{2}\left[  \frac{\lambda}{uz}-\frac{\dot{z}}{z}-\frac
{\mathbf{\dot{J}}}{\mathbf{J}}\right] \\
&  =\frac{1}{2}\left[  \frac{\dot{\omega}}{\omega}-\frac{\dot{z}}{z}%
-\frac{\mathbf{\dot{J}}}{\mathbf{J}}\right]
\end{align*}
where \ $\omega(s,t)$ \ satisfies the ODE%
\[
\frac{d}{dt}\ln(\omega)=\frac{\lambda}{uz}%
\]
which can be solved to give%
\begin{equation}
\omega(s,t)=\frac{\left(  \lambda-s\gamma\right)  u-\lambda}{u\left[
\rho+\gamma(1-\gamma)s\right]  }. \label{w}%
\end{equation}
We conclude that%
\begin{equation}
K(s,t)=K_{0}(s)z^{-\frac{1}{2}}\sqrt{\frac{\omega(s,t)}{\mathbf{J}(s,t)}},
\label{K}%
\end{equation}
with \ $K_{0}(s)$ \ to be determined. As \ $(y,z)\rightarrow(0,\gamma)$
\ (i.e., as \ $t\rightarrow0)$ the Jacobian \ $\mathbf{J}(s,t)\rightarrow0.$
\ Therefore, the asymptotic expansion in (\ref{GRay}) ceases to be valid.

So far we have determined the exponent \ $\psi(s,t)$ \ and the leading
amplitude \ $K(s,t)$ \ except for the constant $\psi_{0}$ in (\ref{psi1}), the
function \ $K_{0}(s)$ in (\ref{K})\ and the power \ $\nu$ \ in (\ref{GRay}).
\ In section 4 we will determine them by matching (\ref{GRay}) to a corner
layer solution valid in a neighborhood of the point \ $(0,\gamma)$.

\subsection{The rays from infinity}

From (\ref{Finf}), we have
\[
F_{k}(\infty)=G(\infty,z)=\exp\left\{  \frac{1}{\varepsilon}\left[
z\ln(\lambda)-\ln(\lambda+1)\right]  \right\}  \binom{\varepsilon^{-1}%
}{z\varepsilon^{-1}},
\]
and by Stirling's formula
\begin{equation}
G(\infty,z)\sim\sqrt{\varepsilon}\kappa(z)\exp\left[  \frac{1}{\varepsilon
}\Phi(z)\right]  ,\quad\varepsilon\rightarrow0,\quad0<z<1 \label{Ginfinity}%
\end{equation}
where%
\begin{equation}
\Phi(z)=-z\ln(z)-(1-z)\ln(1-z)+z\ln(\lambda)-\ln(\lambda+1) \label{Phi}%
\end{equation}%
\begin{equation}
\kappa(z)=\frac{1}{\sqrt{2\pi}}\frac{1}{\sqrt{z(1-z)}}. \label{k}%
\end{equation}
Note that $\Phi(z)$ and $\kappa(z)$ satisfy (\ref{eik}) and (\ref{trans}), respectively.

Denoting the scaled domain by%
\begin{equation}
\mathfrak{D=}[0,\infty)\times\lbrack0,1], \label{D}%
\end{equation}
we must determine what part of $\mathfrak{D}$ the rays from infinity fill.
Defining $p_{\infty}=\frac{\partial\Phi}{\partial y}$ and $u_{\infty}%
=\exp\left(  \frac{\partial\Phi}{\partial z}\right)  ,$ we have%
\begin{equation}
p_{\infty}\equiv0,\quad u_{\infty}=\frac{\lambda\left(  1-z\right)  }{z}.
\label{uinf}%
\end{equation}
(\ref{strip}) and (\ref{eqb}) yield the equations for the rays $y_{\infty
}(t),\ z_{\infty}(t)\ $
\begin{equation}
\dot{y}_{\infty}=z_{\infty}-\gamma,\quad\dot{z}_{\infty}=\left(
1+\lambda\right)  z_{\infty}-\lambda\label{yzinf}%
\end{equation}
or, eliminating $t$ from the system (\ref{yzinf}) and writing $y_{\infty
}(t)=Y_{\infty}(z)$ we get%
\begin{equation}
\frac{dY_{\infty}}{dz}=\frac{z-\gamma}{\left(  1+\lambda\right)  z-\lambda}.
\label{Y(z)}%
\end{equation}
Solving (\ref{yzinf}) subject to the initial condition $Y_{\infty}%
(z_{0})=y_{0},$ where $(y_{0},z_{0})\in\partial\mathfrak{D}$ we get%
\begin{equation}
Y_{\infty}(z)=y_{0}+\frac{z-z_{0}}{\lambda+1}-\frac{\rho}{\left(
\lambda+1\right)  ^{2}}\ln\left[  \frac{(\lambda+1)z-\lambda}{(\lambda
+1)z_{0}-\lambda}\right]  ,\quad z_{0}\neq\frac{\lambda}{\lambda+1}
\label{Yinf}%
\end{equation}
and when $z_{0}=\frac{\lambda}{\lambda+1}$ the ray is a line parallel to the
$y$-axis given by
\begin{equation}
(y_{\infty},\ z_{\infty})\equiv\left(  y,\ \frac{\lambda}{\lambda+1}\right)
,\quad0\leq y. \label{seg}%
\end{equation}
We can interpret (\ref{Yinf}) as a family of rays emanating from the "point"
$\left(  \infty,\frac{\lambda}{\lambda+1}\right)  $ and hitting the domain's
boundary $\partial\mathfrak{D}$ at the point $(y_{0},z_{0}).$ We divide the
rays into five groups (see Figure 3) depending on the location of
$(y_{0},z_{0})$. The quantities inside the logarithms are all greater than one:

\begin{figure}[t]
\begin{center}
\rotatebox{270} {\resizebox{10cm}{!}{\includegraphics{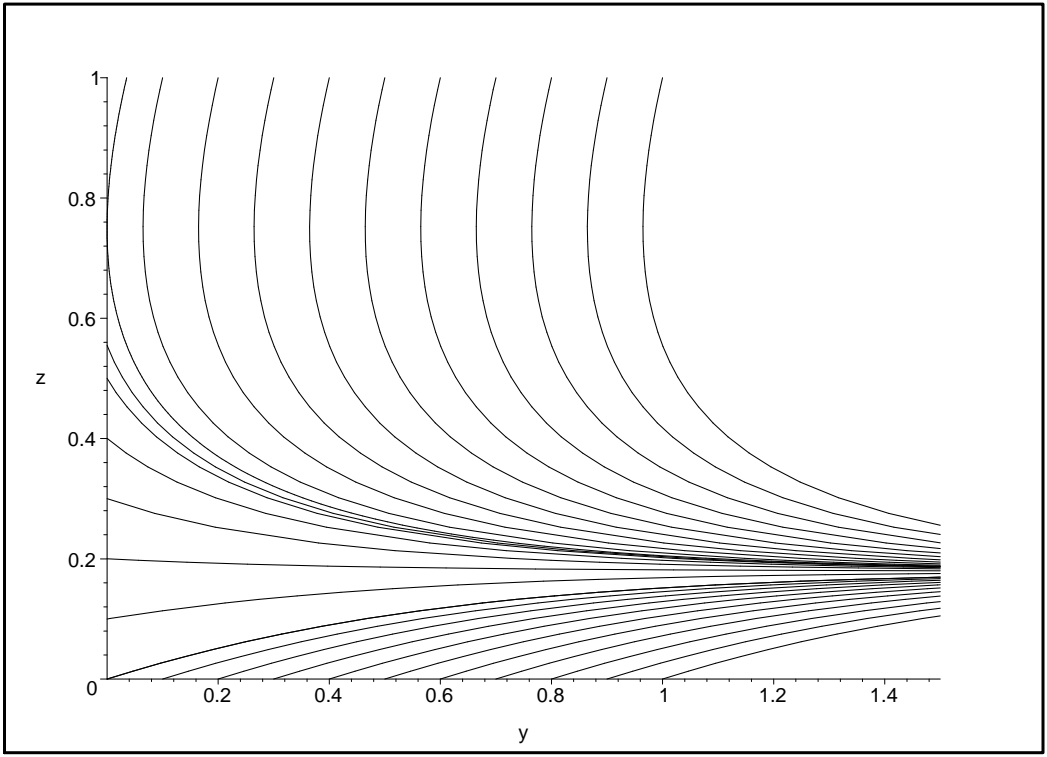}}}
\end{center}
\caption{A sketch of the rays from infinity.}%
\end{figure}

1) $y_{0}\geq0,\quad z_{0}=0.$
\[
Y_{\infty}(z)=y_{0}+\frac{z}{\lambda+1}+\frac{\rho}{\left(  \lambda+1\right)
^{2}}\ln\left[  \frac{\lambda}{\lambda-(\lambda+1)z}\right]  ,\quad0\leq
z<\frac{\lambda}{\lambda+1}.
\]

2) $y_{0}=0,\quad0\leq z_{0}<\frac{\lambda}{\lambda+1}.$
\[
Y_{\infty}(z)=\frac{z-z_{0}}{\lambda+1}+\frac{\rho}{\left(  \lambda+1\right)
^{2}}\ln\left[  \frac{(\lambda+1)z_{0}-\lambda}{(\lambda+1)z-\lambda}\right]
,\quad z_{0}\leq z<\frac{\lambda}{\lambda+1}.
\]
For $z_{0}=\frac{\lambda}{\lambda+1}$ the ray becomes parallel to the $y$-axis
(\ref{seg}).

3) $y_{0}=0,\quad\frac{\lambda}{\lambda+1}<z_{0}<\gamma.$%
\[
Y_{\infty}(z)=\frac{z-z_{0}}{\lambda+1}-\frac{\rho}{\left(  \lambda+1\right)
^{2}}\ln\left[  \frac{(\lambda+1)z_{0}-\lambda}{(\lambda+1)z-\lambda}\right]
,\quad\frac{\lambda}{\lambda+1}<z\leq z_{0}.
\]

4) $y_{0}=0,\quad z_{0}=\gamma.$
\[
Y_{\infty}(z)=\frac{z-\gamma}{\lambda+1}-\frac{\rho}{\left(  \lambda+1\right)
^{2}}\ln\left[  \frac{\rho}{(\lambda+1)z-\lambda}\right]  ,\quad\frac{\lambda
}{\lambda+1}<z\leq1.
\]

This critical ray is tangent to the $z$-axis at $z=\gamma.$ The upper branch
of the ray for $\gamma\leq z\leq1$ is the same as the ray $Y_{0}(z)$ which
emanated from $(0,\gamma)$ and was defined in (\ref{Y0}).

The part of the boundary corresponding to $\left\{  0\right\}  \times
(\gamma,1]\cup\lbrack0,Y_{0}(1))\times\left\{  1\right\}  $ is not reached by
the rays.

5) $y_{0}>Y_{0}(1),\quad z_{0}=1.$%
\[
Y_{\infty}(z)=y_{0}+\frac{z-1}{\lambda+1}-\frac{\rho}{\left(  \lambda
+1\right)  ^{2}}\ln\left[  \frac{1}{(\lambda+1)z-\lambda}\right]  ,\quad
\frac{\lambda}{\lambda+1}<z\leq1.
\]

The rays from infinity fill the region given by%
\begin{equation}
R=\left\{  0\leq y,\quad0\leq z\leq\gamma\right\}  \cup\left\{  Y_{0}(z)\leq
y,\quad\gamma\leq z\leq1\right\}  . \label{R}%
\end{equation}
The complementary region \ $R^{C}$%
\begin{equation}
R^{C}=\left\{  0\leq y<Y_{0}(z),\quad\gamma\leq z\leq1\right\}  , \label{RC}%
\end{equation}
is a \textit{shadow }of the rays from infinity. In $R^{C}$, \ $G$ \ is given
by (\ref{GRay}) as only the rays from $(0,\gamma)$ are present. In the region
\ $R$ \ below$,$ both the rays coming from\ $(0,\gamma)$ and the rays coming
from infinity must be taken into account. We add (\ref{GRay}) and
(\ref{Ginfinity}) to represent \ $G$ \ in the asymptotic form \
\begin{equation}
G(y,z)\sim\sqrt{\varepsilon}\kappa(z)\exp\left[  \frac{1}{\varepsilon}%
\Phi(z)\right]  +\varepsilon^{\nu}\exp\left[  \frac{1}{\varepsilon}%
\Psi(y,z)\right]  \mathbb{K}(y,z),\quad(y,z)\in R. \label{GR}%
\end{equation}
We can show that $\Phi(z)>\Psi(y,z)$ in the interior of $R,$ so that
$G(y,z)\sim G(\infty,z).$ However, in $R$ we can refine (\ref{GR}) to
$G(y,z)-G(\infty,z)\sim\varepsilon^{\nu}\exp\left[  \frac{1}{\varepsilon}%
\Psi(y,z)\right]  \mathbb{K}(y,z).$

\section{The corner layer at \ $(0,\gamma)$}

We determine the constant $\nu$ in (\ref{GRay}) and the function $K_{0}(s)$ in
(\ref{K}) by considering carefully the region where the rays from $(0,\gamma)$
enter the domain $\mathfrak{D}$ and using asymptotic matching. We introduce
the stretched variables $\chi,\ l,$ the function $G_{1}(l,\chi)$ and the
parameter $\alpha$ defined by%
\begin{align}
\chi &  =\frac{x}{\varepsilon},\quad\chi\geq0,\quad F_{k}(x)=G_{1}%
(k-c+\alpha,\frac{x}{\varepsilon})=G_{1}(l,\chi)\nonumber\\
l  &  =k-c+\alpha,\quad-\infty<l<\infty\label{alpha}\\
\alpha &  =c-\left\lfloor c\right\rfloor ,\quad0<\alpha<1.\nonumber
\end{align}
Note that $\alpha$ is the fractional part of $c.$ Use of (\ref{alpha}) in
(\ref{diffeq}) yields the equation%
\begin{gather*}
(l-\alpha)\frac{dG_{1}}{d\chi}\varepsilon^{-1}=\left[  \gamma G_{1}\left(
l+1,\chi\right)  +\lambda(1-\gamma)G_{1}\left(  l-1,\chi\right)  -\phi
G_{1}\left(  l,\chi\right)  \right]  \varepsilon^{-1}\\
+\lambda(\alpha+1-l)G_{1}\left(  l-1,\chi\right)  +(l-\alpha+1)G_{1}\left(
l+1,\chi\right)  +\left(  \lambda-1\right)  (l-\alpha)G_{1}\left(
l,\chi\right)
\end{gather*}
or, to leading order in $\varepsilon,$%
\begin{equation}
(l-\alpha)\frac{F_{l}^{(1)}}{d\chi}=\gamma F_{l+1}^{(1)}+\lambda
(1-\gamma)F_{l-1}^{(1)}-\phi F_{l}^{(1)} \label{diffcorner}%
\end{equation}
with $G_{1}\left(  l,\chi\right)  \sim F_{l}^{(1)}\left(  \chi\right)  $ as
$\varepsilon\rightarrow0$ and $\phi=\gamma+\lambda-\gamma\lambda$.

Also (\ref{FBC}) gives the boundary condition%
\begin{equation}
F_{l}^{(1)}(0)=0,\quad l\geq1 \label{BCcorner}%
\end{equation}
and (\ref{Ginfinity}) implies that $F_{k}(\infty)\sim F_{l}^{(1)}(\infty)$%
\begin{equation}
F_{l}^{(1)}(\infty)=\sqrt{\varepsilon}\left(  u_{0}\right)  ^{l-\alpha}%
\kappa(\gamma)\exp\left[  \frac{1}{\varepsilon}\Phi(\gamma)\right]
\label{infcorner}%
\end{equation}
with $u_{0}=\frac{\lambda}{\gamma}\left(  1-\gamma\right)  .$ For a fixed $l$
and $\chi\rightarrow\infty,$ we approach the interior of $R,$ where (\ref{GR})
applies. Note that $F_{l}^{(1)}(\infty)$ is the same as $\sqrt{\varepsilon
}\kappa(z)\exp\left[  \frac{1}{\varepsilon}\Phi(z)\right]  $ expanded for
$z\rightarrow\gamma.$ Thus, (\ref{infcorner}) is the asymptotic matching
condition between the corner layer and the solution in $R.$ We shall examine
the matching to $R^{C}$ later.

Equation (\ref{diffcorner}) admits the separable solutions%
\begin{equation}
F_{l}^{(1)}(\chi)=e^{\theta\chi}h_{l}(\theta) \label{Fl}%
\end{equation}
if \ $h_{l}(\theta)$\ satisfies the difference equation%
\[
\gamma h_{l+1}+\lambda(1-\gamma)h_{l-1}=\left[  (l-\alpha)\theta+\phi\right]
h_{l}.
\]
Setting \ $h_{l}(\theta)=\left(  u_{0}\right)  ^{\frac{l}{2}}H_{l}(\theta)$ we
see that \
\begin{equation}
H_{l+1}+H_{l-1}=\frac{2}{\beta}\left[  (l-\alpha)\theta+\phi\right]  H_{l}
\label{H}%
\end{equation}
with%
\begin{equation}
\beta=2\sqrt{\lambda\gamma\left(  1-\gamma\right)  }. \label{beta}%
\end{equation}
The only solutions to (\ref{H}) which have acceptable behavior as
$l\rightarrow\infty$ are of the form%
\[
H_{l}(\theta)=J_{l-\alpha+\frac{\phi}{\theta}}\left(  \frac{\beta}{\theta
}\right)
\]
where $J$ is the Bessel function. If (\ref{Fl}) is not to grow as
\ $\chi\rightarrow\infty$, we need $\theta\leq0.$ But except when $v$ is an
integer, the Bessel function $J_{v}(x)$ is complex for negative argument.
Therefore, we need%
\[
\frac{\phi}{\theta}-\alpha=-1,-2,\ldots
\]
or
\begin{equation}
\theta_{j}=-\frac{\phi}{j+1-\alpha}<0,\quad j\geq0. \label{thetaj}%
\end{equation}
It follows that the general solution to (\ref{diffcorner}) takes the form%
\[
F_{l}^{(1)}(\chi)=F_{l}^{(1)}(\infty)+\left(  u_{0}\right)  ^{\frac{l}{2}}%
{\displaystyle\sum\limits_{j\geq0}}
a_{j}e^{\theta_{j}\chi}J_{l-\alpha+\frac{\phi}{\theta_{j}}}\left(  \frac
{\beta}{\theta_{j}}\right)
\]
or%
\begin{align}
F_{l}^{(1)}(\chi)  &  =\sqrt{\varepsilon}\left(  u_{0}\right)  ^{l-\alpha
}\kappa(\gamma)\exp\left[  \frac{1}{\varepsilon}\Phi(\gamma)\right]
\label{Fcorner}\\
&  +\left(  \sqrt{u_{0}}\right)  ^{l}%
{\displaystyle\sum\limits_{j\geq0}}
a_{j}\exp\left(  -\frac{\phi}{j+1-\alpha}\chi\right)  J_{l-1-j}\left[
-\frac{\beta}{\phi}(j+1-\alpha)\right] \nonumber
\end{align}
where the coefficients \ $a_{j}$ \ in the above (spectral) representation
remain to be determined.

Taking the Laplace transform
\[
\hat{F}_{l}^{(1)}(\vartheta)=%
{\displaystyle\int\limits_{0}^{\infty}}
e^{-\vartheta\chi}F_{l}^{(1)}(\chi)d\chi
\]
of (\ref{Fcorner}) we obtain%
\begin{align}
\hat{F}_{l}^{(1)}(\vartheta)  &  =\sqrt{\varepsilon}\left(  u_{0}\right)
^{l-\alpha}\kappa(\gamma)\exp\left[  \frac{1}{\varepsilon}\Phi(\gamma)\right]
\frac{1}{\vartheta}\label{lap}\\
&  +\left(  \sqrt{u_{0}}\right)  ^{l}%
{\displaystyle\sum\limits_{j\geq0}}
a_{j}\frac{1}{\vartheta+\frac{\phi}{j+1-\alpha}}J_{j+1-l}\left[  \frac{\beta
}{\phi}(j+1-\alpha)\right]  .\nonumber
\end{align}
Thus, the only singularities of \ $\hat{F}_{l}^{(1)}(\vartheta)$ \ are simple
poles at \ $\vartheta=0$ and \ $\vartheta=\theta_{j},\quad\ j\geq0.$ \ It is
well known that the gamma function \ $\Gamma(z)$ has simple poles at
\ $z=0,-1,-2,\ldots$. Hence, we shall represent \ $\hat{F}_{l}(\vartheta)$ as%
\begin{equation}
\hat{F}_{l}^{(1)}(\vartheta)=\left(  \sqrt{u_{0}}\right)  ^{l}\frac
{1}{\vartheta}\Gamma\left(  \frac{\phi}{\vartheta}+1-\alpha\right)
J_{l-\alpha+\frac{\phi}{\vartheta}}\left(  \frac{\beta}{\vartheta}\right)
f(\vartheta) \label{laplace1}%
\end{equation}
where \ $f(\vartheta)$ is chosen such that $\Gamma\left(  \frac{\phi
}{\vartheta}+1-\alpha\right)  J_{l-\alpha+\frac{\phi}{\vartheta}}\left(
\frac{\beta}{\vartheta}\right)  f(\vartheta)$ is analytic for
\ $\operatorname{Re}(\vartheta)>-\frac{\phi}{1-\alpha}$. Taking the Laplace
transform in (\ref{Fcorner}) we get the equation%
\[
(l-\alpha)\vartheta\hat{F}_{l}^{(1)}=\gamma\hat{F}_{l+1}^{(1)}+\lambda
(1-\gamma)\hat{F}_{l-1}^{(1)}-\phi\hat{F}_{l}^{(1)},\ l\geq1
\]
which is satisfied by (\ref{laplace1}). By the inversion formula of the
Laplace transform we have%
\[
F_{l}^{(1)}(\chi)=\left(  \sqrt{u_{0}}\right)  ^{l}\frac{1}{2\pi i}%
{\displaystyle\int\limits_{\mathrm{Br}}}
e^{\chi\vartheta}\frac{1}{\vartheta}\Gamma\left(  \frac{\phi}{\vartheta
}+1-\alpha\right)  J_{l-\alpha+\frac{\phi}{\vartheta}}\left(  \frac{\beta
}{\vartheta}\right)  f(\vartheta)d\vartheta
\]
where $\mathrm{Br}$ is a vertical contour in the complex \ $\vartheta-$plane
on which \ $\operatorname{Re}(\vartheta)>0.$

Since the residue of \ $\hat{F}_{l}^{(1)}(\vartheta)$ at \ $\vartheta=0$
\ corresponds to \ $F_{l}^{(1)}(\infty)$, we must have%
\[
\left(  \sqrt{u_{0}}\right)  ^{l}\Gamma\left(  \frac{\phi}{\vartheta}%
+1-\alpha\right)  J_{l-\alpha+\frac{\phi}{\vartheta}}\left(  \frac{\beta
}{\vartheta}\right)  f(\vartheta)\rightarrow\sqrt{\varepsilon}\left(
u_{0}\right)  ^{l-\alpha}\kappa(\gamma)\exp\left[  \frac{1}{\varepsilon}%
\Phi(\gamma)\right]
\]
as \ $\vartheta\rightarrow0.$ In that limit we see that%
\[
\Gamma\left(  \frac{\phi}{\vartheta}+1-\alpha\right)  J_{l-\alpha+\frac{\phi
}{\vartheta}}\left(  \frac{\beta}{\vartheta}\right)  \sim\sqrt{\frac
{\vartheta}{\rho}}\left(  \sqrt{u_{0}}\right)  ^{\frac{\phi}{\vartheta
}+l-\alpha}e^{\frac{\rho-\phi}{\vartheta}}\left(  \frac{\phi}{\vartheta
}\right)  ^{\frac{\phi}{\vartheta}+\frac{1}{2}-\alpha},\quad\vartheta
\rightarrow0.
\]
Therefore, we write%
\begin{equation}
f(\vartheta)=\sqrt{\varepsilon}\left(  \sqrt{u_{0}}\right)  ^{-\alpha}%
\sqrt{\frac{\rho}{\phi}}\kappa(\gamma)\exp\left[  \frac{1}{\varepsilon}%
\Phi(\gamma)+\Upsilon(\vartheta)\right]  \tilde{f}(\vartheta) \label{f}%
\end{equation}
where%
\begin{equation}
\Upsilon(\vartheta)=\left(  \frac{\phi}{\vartheta}-\alpha\right)  \ln\left(
\frac{\vartheta}{\phi}\right)  +\frac{2\lambda\left(  1-\gamma\right)
}{\vartheta}-\frac{\phi}{2\vartheta}\ln\left(  u_{0}\right)  \label{Lambda}%
\end{equation}
and\ \ $\tilde{f}(\vartheta)$ is entire, with \ $\tilde{f}(0)=1.$

By combining the preceding results we have%
\begin{align}
F_{l}^{(1)}(\chi)  &  =\sqrt{\varepsilon}\sqrt{\frac{\rho}{\phi}}\kappa
(\gamma)\left(  \sqrt{u_{0}}\right)  ^{l-\alpha}\exp\left[  \frac
{1}{\varepsilon}\Phi(\gamma)\right] \label{Fcorner1}\\
&  \times\frac{1}{2\pi i}%
{\displaystyle\int\limits_{Br}}
e^{\chi\vartheta}\frac{1}{\vartheta}\Gamma\left(  \frac{\phi}{\vartheta
}+1-\alpha\right)  J_{l-\alpha+\frac{\phi}{\vartheta}}\left(  \frac{\beta
}{\vartheta}\right)  \exp\left[  \Upsilon(\vartheta)\right]  \tilde
{f}(\vartheta)d\vartheta.\nonumber
\end{align}
The boundary condition (\ref{BCcorner}) implies that
\[
\underset{\vartheta\rightarrow\infty}{\lim}\left[  \vartheta\hat{F}_{l}%
^{(1)}(\vartheta)\right]  =0,\quad l\geq1.
\]
From (\ref{Fcorner1}) we have%
\[
\frac{1}{\vartheta}\Gamma\left(  \frac{\phi}{\vartheta}+1-\alpha\right)
J_{l-\alpha+\frac{\phi}{\vartheta}}\left(  \frac{\beta}{\vartheta}\right)
\exp\left[  \Upsilon(\vartheta)\right]  \sim\left(  \frac{1}{\vartheta
}\right)  ^{l+1}\frac{\Gamma\left(  1-\alpha\right)  }{\Gamma\left(
l-\alpha+1\right)  }\left(  \frac{\beta}{2}\right)  ^{l-\alpha}\phi^{\alpha
},\quad\vartheta\rightarrow\infty.
\]
Fixing $l=1$, we get \ $\tilde{f}(\vartheta)=o\left(  \vartheta\right)
,\quad\vartheta\rightarrow\infty$ and by Liouville's theorem this forces
\ $\tilde{f}(\vartheta)$ to be a constant. Since \ $\tilde{f}(0)=1,$ we have
$\tilde{f}(\vartheta)\equiv1.$ Thus, (\ref{Fcorner1}) becomes%
\begin{align}
F_{l}^{(1)}(\chi)  &  =\sqrt{\varepsilon}\sqrt{\frac{\rho}{\phi}}\kappa
(\gamma)\left(  \sqrt{u_{0}}\right)  ^{l-\alpha}\exp\left[  \frac
{1}{\varepsilon}\Phi(\gamma)\right] \label{Fcorner2}\\
&  \times\frac{1}{2\pi i}%
{\displaystyle\int\limits_{Br}}
e^{\chi\vartheta}\frac{1}{\vartheta}\Gamma\left(  \frac{\phi}{\vartheta
}+1-\alpha\right)  J_{l-\alpha+\frac{\phi}{\vartheta}}\left(  \frac{\beta
}{\vartheta}\right)  \exp\left[  \Upsilon(\vartheta)\right]  d\vartheta
.\nonumber
\end{align}
The coefficients \ $a_{j}$ in the spectral expansion (\ref{Fcorner}) are
determined from (\ref{Fcorner2}) by applying the residue theorem. Noting that%
\[
\operatorname*{Res}\left[  \Gamma\left(  \frac{\phi}{\vartheta}+1-\alpha
\right)  ,\ \vartheta=\theta_{j}\right]  =-\frac{\phi}{\left(  j+1-\alpha
\right)  ^{2}}\frac{\left(  -1\right)  ^{j}}{j!}%
\]
we obtain%
\begin{equation}
\ a_{j}=\sqrt{\varepsilon}\sqrt{\frac{\rho}{\phi}}\kappa(\gamma)\left(
\sqrt{u_{0}}\right)  ^{-\alpha}\exp\left[  \frac{1}{\varepsilon}\Phi
(\gamma)+\Upsilon(\theta_{j})\right]  \frac{\left(  -1\right)  ^{j}}{j!}%
\frac{1}{\left(  j+1-\alpha\right)  },\quad j\geq0. \label{aj}%
\end{equation}
This completes the determination of the spectral and integral representations
of $F_{l}^{(1)}(\chi)$ and hence the leading term for $F_{k}(x)$ in the corner region.

\subsection{\ \ \ Matching the corner and $R^{C}$ regions}

In this section we shall determine the constant $\psi_{0}$ in (\ref{psi1}),
the function \ $K_{0}(s)$ in (\ref{K})\ and the power \ $\nu$ \ in
(\ref{GRay}). From (\ref{Y0}) we have%
\[
Y_{0}(z)\sim\frac{1}{2\rho}(z-\gamma)^{2},\quad z\rightarrow\gamma.
\]
If we introduce the new variable \ $\Omega$ defined by
\begin{equation}
\Omega=\frac{2\rho y}{(z-\gamma)^{2}}=\frac{2\rho\chi}{(l-\alpha)^{2}}%
,\quad\Omega=O(1) \label{Omega}%
\end{equation}
we see that the limit \ $\chi,\ l\rightarrow\infty,\ \Omega>1,$ corresponds to
the matching between the corner and $R^{C}$ regions.

We set%
\[
\vartheta=\varepsilon\Theta,\quad\eta=(l-\alpha)\vartheta+\phi=(z-\gamma
)\Theta+\phi,\quad\eta,\Theta=O(1),\quad\eta,\Theta>0
\]
use this in \ $J_{l-\alpha+\frac{\phi}{\vartheta}}\left(  \frac{\beta
}{\vartheta}\right)  $ and let \ $\varepsilon\rightarrow0,$ thus obtaining%
\[
J_{l-\alpha+\frac{\phi}{\vartheta}}\left(  \frac{\beta}{\vartheta}\right)
=J_{\upsilon}\left(  \frac{\beta}{\eta}\upsilon\right)  \sim\frac{\exp\left\{
\upsilon\left[  \sqrt{1-\left(  \beta/\eta\right)  ^{2}}+\ln\left(
\frac{\beta/\eta}{1+\sqrt{1-\left(  \beta/\eta\right)  ^{2}}}\right)  \right]
\right\}  }{\sqrt{2\pi\upsilon\sqrt{1-\left(  \beta/\eta\right)  ^{2}}}}%
\]%
\[
=\frac{\sqrt{\varepsilon\Theta}}{\sqrt{2\pi}\sqrt{p(\eta)}}\exp\left\{
\frac{1}{\varepsilon\Theta}\left[  p(\eta)-\eta\ln\left(  \frac{\eta+p(\eta
)}{\beta}\right)  \right]  \right\}
\]
with%
\begin{equation}
\upsilon=\frac{\eta}{\varepsilon\Theta},\quad p(\eta)=\sqrt{\eta^{2}-\beta
^{2}},\quad p(\phi)=\rho. \label{P}%
\end{equation}
Use of Stirling's formula gives%
\[
\Gamma\left(  \frac{\phi}{\varepsilon\Theta}+1-\alpha\right)  \sim\sqrt{2\pi
}\exp\left\{  \frac{\phi}{\varepsilon\Theta}\left[  \ln\left(  \frac{\phi
}{\varepsilon\Theta}\right)  -1\right]  \right\}  \left(  \frac{\phi
}{\varepsilon\Theta}\right)  ^{\frac{1}{2}-\alpha}%
\]
and from (\ref{Lambda}) we have%
\[
\exp\left[  \Upsilon(\varepsilon\Theta)\right]  =\exp\left\{  \frac
{1}{\varepsilon\Theta}\left[  \phi\ln\left(  \frac{\varepsilon\Theta}{\phi
}\right)  +2\lambda\left(  1-\gamma\right)  -\frac{\phi}{2}\ln\left(
u_{0}\right)  \right]  \right\}  \left(  \frac{\varepsilon\Theta}{\phi
}\right)  ^{-\alpha}.
\]
Therefore,%
\begin{gather}
J_{\upsilon}\left(  \frac{\beta}{\eta}\upsilon\right)  \Gamma\left(
\frac{\phi}{\varepsilon\Theta}+1-\alpha\right)  \exp\left[  \Upsilon
(\varepsilon\Theta)\right]  \sim\label{asymp1}\\
\frac{\sqrt{\phi}}{\sqrt{p(\eta)}}\exp\left\{  \frac{1}{\varepsilon\Theta
}\left[  p(\eta)-\eta\ln\left(  \frac{\eta+p(\eta)}{\beta}\right)  -\rho
-\frac{\phi}{2}\ln\left(  u_{0}\right)  \right]  \right\}  .\nonumber
\end{gather}

Using (\ref{asymp1}) in (\ref{Fcorner2}) yields, in terms of $z$ and $\Omega,$%
\begin{align}
F_{l}^{(1)}(\chi)  &  \sim\sqrt{\varepsilon}\sqrt{\rho}\kappa(\gamma
)\exp\left\{  \frac{1}{\varepsilon}\left[  \Phi(\gamma)+\frac{z-\gamma}{2}%
\ln\left(  u_{0}\right)  \right]  \right\} \label{asymp2}\\
&  \times\frac{1}{2\pi i}%
{\displaystyle\int\limits_{\mathrm{Br}^{\prime}}}
\frac{1}{\eta-\phi}\frac{1}{\sqrt{p(\eta)}}\exp\left[  \frac{1}{\varepsilon
}\left(  z-\gamma\right)  g(\eta)\right]  d\eta\nonumber
\end{align}
where%
\begin{equation}
g(\eta)=\frac{\left(  \eta-\phi\right)  \Omega}{2\rho}+\frac{1}{\eta-\phi
}\left[  p(\eta)-\eta\ln\left(  \frac{\eta+p(\eta)}{\beta}\right)  -\rho
-\frac{\phi}{2}\ln\left(  u_{0}\right)  \right]  \label{g}%
\end{equation}
and $\mathrm{Br}^{\prime}$ is a vertical contour in the complex plane with
$\operatorname{Re}(\eta)>\phi.$ For \ $\varepsilon\rightarrow0$ \ with
$\Omega$ \ fixed, we can evaluate (\ref{asymp2}) by the saddle point method to
get%
\begin{gather}
F_{l}^{(1)}(\chi)\sim\sqrt{\varepsilon}\sqrt{\rho}\kappa(\gamma)\exp\left\{
\frac{1}{\varepsilon}\left[  \Phi(\gamma)+\frac{z-\gamma}{2}\ln\left(
u_{0}\right)  \right]  \right\} \label{saddle}\\
\times\frac{\sqrt{\varepsilon}}{\sqrt{2\pi}}\frac{1}{\sqrt{z-\gamma}}\frac
{1}{\eta^{\ast}-\phi}\frac{1}{\sqrt{p(\eta^{\ast})}}\exp\left[  \frac
{1}{\varepsilon}\left(  z-\gamma\right)  g(\eta^{\ast})\right]  \frac{1}%
{\sqrt{g^{\prime\prime}(\eta^{\ast})}}\nonumber
\end{gather}
where the saddle point \ $\eta^{\ast}\left(  \Omega\right)  $ \ is defined by
\ $g^{\prime}(\eta^{\ast})=0$. Note that $\eta^{\ast}\left(  \Omega\right)
>\phi$ for $\Omega>1,$ i.e., the saddle point $\eta^{\ast}\left(
\Omega\right)  $ lies to the right of the pole at $\eta=\phi$ and the
integrand is analytic for $\operatorname{Re}(\eta)>\phi$.

Taking derivatives in (\ref{g}) we find that%

\begin{equation}
g^{\prime}(\eta)=\frac{\Omega}{2\rho}+\frac{1}{\left(  \eta-\phi\right)  ^{2}%
}\left[  \rho+\frac{\phi}{2}\ln\left(  u_{0}\right)  +\phi\ln\left(
\frac{\eta+p(\eta)}{\beta}\right)  -p(\eta)\right]  \label{dg}%
\end{equation}%
\begin{equation}
g^{\prime\prime}(\eta)=\frac{1}{\left(  \eta-\phi\right)  ^{3}}\left[
\frac{\eta^{2}+2\eta\phi-\phi^{2}-2\beta^{2}}{p(\eta)}-2\phi\ln\left(
\frac{\eta+p(\eta)}{\beta}\right)  -2\rho-\phi\ln\left(  u_{0}\right)
\right]  . \label{d2g}%
\end{equation}

From (\ref{dg}) we observe that \ $g^{\prime}(\eta^{\ast})=0$ \ if
\ $\eta^{\ast}=\phi$ and $\Omega=1,$ which implies that \ $\eta^{\ast}\left(
1\right)  =\phi.$ To determine $\eta^{\ast}$ for $\Omega\sim1,$ we use
(\ref{dg}) and an expansion of the form%
\begin{equation}
\eta^{\ast}\left(  \Omega\right)  \sim\phi+a_{1}\left(  \Omega-1\right)
+a_{2}\left(  \Omega-1\right)  ^{2}+a_{3}\left(  \Omega-1\right)  ^{3}+\cdots.
\label{etastar}%
\end{equation}
Using (\ref{etastar}) in (\ref{dg}) and expanding the latter in powers of
\ $\Omega-1,$ we find that
\[
a_{1}=-\frac{3\rho^{2}}{2\phi},\quad a_{2}=-\frac{27\rho^{2}}{32\phi^{3}%
}\left(  \rho^{2}-3\phi^{2}\right)
\]
and
\begin{align*}
g(\eta^{\ast})  &  \sim\frac{1}{2}\ln\left(  u_{0}\right)  -\frac{3\rho}%
{8\phi}\left(  \Omega-1\right)  ^{2},\quad g^{\prime\prime}(\eta^{\ast}%
)\sim\frac{\phi}{3\rho^{3}},\\
\quad\frac{1}{\eta^{\ast}-\phi}  &  \sim-\frac{2\phi}{3\rho^{2}}\left(
\Omega-1\right)  ^{-1},\quad\frac{1}{\sqrt{p(\eta^{\ast})}}\sim\frac{1}%
{\sqrt{\rho}}%
\end{align*}
from which we conclude that as $(y,z)\rightarrow(0,\gamma)$%
\begin{equation}
F_{l}^{(1)}(\chi)\sim-\frac{1}{\sqrt{3}}\frac{\varepsilon}{\pi}\frac
{\sqrt{\phi}}{\sqrt{\rho}}\frac{1}{\sqrt{\gamma\left(  1-\gamma\right)  }%
}\frac{1}{\sqrt{z-\gamma}}\left(  \Omega-1\right)  ^{-1}\exp\left\{  \frac
{1}{\varepsilon}\left[  \Phi(\gamma)+\left(  z-\gamma\right)  \ln\left(
u_{0}\right)  \right]  \right\}  . \label{fsaddle}%
\end{equation}

We next evaluate $K$ and $\psi$ in (\ref{GRay}) near the corner $(0,\gamma
).$Using $T_{+}(s,z),$ (cf. (\ref{T})) in (\ref{y(s,t)}) and expanding for
small \ $s$, we obtain%
\begin{equation}
y\sim Y_{0}(z)-Y_{2}(z)s \label{y11}%
\end{equation}
with%
\begin{gather}
Y_{2}(z)=\frac{2\zeta}{\left(  \lambda+1\right)  ^{4}}\ln\left(
\frac{z+z\lambda-\lambda}{\rho}\right) \label{Y2}\\
-\frac{z-\gamma}{\left(  \lambda+1\right)  \left(  \lambda z+z-\lambda\right)
^{2}}\left[  \frac{2\zeta\rho}{\left(  \lambda+1\right)  ^{2}}+\frac{3\zeta
}{\left(  \lambda+1\right)  }\left(  z-\gamma\right)  +\left(  \lambda
-1\right)  \left(  z-\gamma\right)  ^{2}\right]  ,\nonumber
\end{gather}
and%
\begin{equation}
\zeta=2\lambda-\gamma+(\gamma-1)\lambda^{2}. \label{zeta}%
\end{equation}
When \ $z$ \ is close to \ $\gamma,$ (\ref{y11}) yields%
\[
y\sim\frac{1}{2\rho}\left(  z-\gamma\right)  ^{2}-\frac{\left(  \lambda
+1\right)  \rho+\phi s}{3\rho^{3}}\left(  z-\gamma\right)  ^{3},\quad
z\rightarrow\gamma
\]
or, using (\ref{Omega}),
\begin{equation}
s\sim-\frac{3}{2}\frac{\rho^{2}}{\phi}\frac{\Omega-1}{z-\gamma}-\left(
\lambda+1\right)  \frac{\rho}{\phi},\quad z\rightarrow\gamma. \label{s}%
\end{equation}
Also, from (\ref{T}) we have%
\begin{equation}
t\sim\frac{z-\gamma}{\rho},\quad z\rightarrow\gamma. \label{t1}%
\end{equation}

We expand (\ref{psi1}) for small \ $t$%
\[
\psi(s,t)\sim\psi_{0}+\ln\left(  u_{0}\right)  \rho t,\quad t\rightarrow0
\]
which taking (\ref{t1}) into account gives%
\[
\Psi(y,z)\sim\psi_{0}+\left(  z-\gamma\right)  \ln\left(  u_{0}\right)  ,\quad
z\rightarrow\gamma
\]
in agreement with (\ref{fsaddle}) if
\begin{equation}
\psi_{0}=\Phi(\gamma). \label{psi0}%
\end{equation}

From (\ref{K}) we obtain%
\[
K(s,t)\sim K_{0}(s)\frac{1}{\sqrt{\gamma}}\sqrt{\frac{3}{\left(
1-\gamma\right)  \phi\rho t^{3}}},\quad t\rightarrow0
\]
or, using (\ref{s}) and (\ref{t1}) in the above,%
\begin{equation}
\mathbb{K}(y,z)\sim K_{0}\left[  -\frac{3}{2}\frac{\rho^{2}}{\phi}\frac
{\Omega-1}{z-\gamma}\right]  \frac{\rho}{\sqrt{\gamma}}\sqrt{\frac{3}{\left(
1-\gamma\right)  \phi\left(  z-\gamma\right)  ^{3}}},\quad z\rightarrow\gamma.
\label{K1}%
\end{equation}
Matching the algebraic factors in (\ref{fsaddle}) and (\ref{K1}) yields%
\[
K_{0}\left[  -\frac{3}{2}\frac{\rho^{2}}{\phi}\frac{\Omega-1}{z-\gamma
}\right]  =-\frac{1}{3\pi}\frac{\phi}{\rho^{\frac{3}{2}}}\frac{z-\gamma
}{\Omega-1}%
\]
which implies that%
\begin{equation}
K_{0}(s)=\frac{\sqrt{\rho}}{2\pi s}. \label{K0}%
\end{equation}
Finally, the exponent of \ $\varepsilon$ \ in (\ref{GRay}) is determined from
(\ref{fsaddle}) and turns out to be \ \ $\nu=1.$ This completes the
determination of the asymptotic solution corresponding to rays from the point
\ $(0,\gamma).$ To summarize, we have established the following.

\begin{result}
The solution of (\ref{eqG}) is asymptotically given by%
\begin{equation}
G(y,z)\sim\varepsilon\exp\left[  \frac{1}{\varepsilon}\Psi(y,z)\right]
\mathbb{K}(y,z)\text{ \ in }R^{C} \label{G4}%
\end{equation}%
\begin{equation}
G(\infty,z)-G(y,z)\sim-\varepsilon\exp\left[  \frac{1}{\varepsilon}%
\Psi(y,z)\right]  \mathbb{K}(y,z)\text{ \ in }R \label{G3}%
\end{equation}
with%
\[
G(\infty,z)=\exp\left\{  \frac{1}{\varepsilon}\left[  z\ln(\lambda
)-\ln(\lambda+1)\right]  \right\}  \binom{\varepsilon^{-1}}{z\varepsilon^{-1}}%
\]%
\[
R=\left\{  0\leq y,\quad0\leq z\leq\gamma\right\}  \cup\left\{  Y_{0}(z)\leq
y,\quad\gamma\leq z\leq1\right\}
\]%
\[
Y_{0}(z)=\frac{z-\gamma}{\lambda+1}-\frac{\rho}{\left(  \lambda+1\right)
^{2}}\ln\left[  \frac{\rho}{(\lambda+1)z-\lambda}\right]  ,\quad\gamma\leq
z\leq1
\]%
\[
\rho=\gamma\lambda+\gamma-\lambda,\quad u_{0}=\frac{\lambda}{\gamma}\left(
1-\gamma\right)
\]%
\begin{align}
\Psi(y,z)  &  =\psi(s,t)=\ln\left[  \frac{\lambda}{\gamma\left(
\lambda+1\right)  }\right]  +sy+(z-1)\ln(u)+\frac{(\lambda-\gamma
s)r_{2}-\lambda}{r_{2}\Delta}\ln\left(  \frac{u-r_{2}}{u_{0}-r_{2}}\right)
\label{psi3}\\
&  -\frac{(\lambda-\gamma s)r_{1}-\lambda}{r_{1}\Delta}\ln\left(
\frac{u-r_{1}}{u_{0}-r_{1}}\right) \nonumber
\end{align}%
\begin{equation}
\mathbb{K}(y,z)=K(s,t)=\frac{\sqrt{\rho}}{2\pi s}z^{-\frac{1}{2}}\sqrt
{\frac{\omega(s,t)}{\mathbf{J}(s,t)}} \label{K2}%
\end{equation}%
\[
r_{1,2}(s)=\frac{1}{2}\left(  s+1-\lambda\pm\Delta\right)  ,\quad
\Delta(s)=\sqrt{\left(  \lambda-s-1\right)  ^{2}+4\lambda},
\]
$(y,z)$ is related to $(s,t)$ by (\ref{y(s,t)}), (\ref{z(s,t)}) and
$u(s,t),\ \omega(s,t),\ \mathbf{J}(s,t)$ are defined in (\ref{u(s,t)}),
(\ref{w}) and (\ref{J}) respectively. We note that $s<0$ in $R$ so that the
right side of (\ref{G3}) is positive. This gives the leading term for the
probability%
\[
\Pr\left[  X(\infty)>x=\frac{y}{\varepsilon},\quad Z(\infty)=k=\frac
{z}{\varepsilon}\right]
\]
that the buffer exceeds $x=Ny.$

In the corner range where (\ref{alpha}) applies, the leading term is given by
(\ref{Fcorner2}), or (\ref{Fcorner}) and (\ref{aj}).
\end{result}

\section{Transition layer}

We shall find a transition layer solution near the curve \ $y=Y_{0}(z)$
defined by (\ref{Y0}) which separates \ $R$ and $R^{C}.$ On this curve
\ $s=0,$ hence (\ref{K2}) is not valid because \ $\mathbb{K}(y,z)$ \ is
infinite there.

We introduce the new function \ $L_{k}(x)$ defined by%
\[
F_{k}(x)=F_{k}(\infty)L_{k}(x).
\]
Then (\ref{diffeq}) yields for \ $L_{k}(x)$ the equation%
\[
(k-c)L_{k}^{\prime}=\left[  \lambda\left(  k-N\right)  -k\right]
L_{k}+\lambda\left(  N-k\right)  L_{k+1}+kL_{k-1}%
\]
and matching to region $R$ forces%
\begin{equation}
L_{k}(\infty)=1. \label{Linf}%
\end{equation}
In terms of the variables \ $y=\varepsilon x,\quad z=\varepsilon k,$ the
function $L^{(1)}(y,z)=$ $L_{k}(x)$ and the parameters \ $\gamma=\varepsilon
c,\quad\varepsilon=\frac{1}{N},$ we get%
\[
(z-\gamma)\frac{\partial L^{(1)}}{\partial y}=(z+\lambda z-\lambda
)\frac{\partial L^{(1)}}{\partial z}+\frac{\varepsilon}{2}(z-\lambda
z+\lambda)\frac{\partial^{2}L^{(1)}}{\partial z^{2}}+O\left(  \frac
{\partial^{3}L^{(1)}}{\partial z^{3}}\varepsilon^{2}\right)  .
\]
Introducing the stretched variable \ $\Lambda$ \ defined by%
\begin{equation}
y=Y_{0}(z)+\sqrt{\varepsilon}\Lambda\label{lambda}%
\end{equation}
and the function $L^{(2)}(\Lambda,z)=L^{(1)}(y,z),$ we obtain for
\ $L^{(2)}(\Lambda,z)$ \ to leading order the diffusion equation%
\begin{equation}
(\lambda-z-\lambda z)\frac{\partial L^{(2)}}{\partial z}-\frac{1}{2}%
\frac{\lambda z-z-\lambda}{\left(  \lambda z+z-\lambda\right)  ^{2}}\left(
z-\gamma\right)  ^{2}\frac{\partial^{2}L^{(2)}}{\partial\Lambda^{2}}=0.
\label{L}%
\end{equation}
To solve (\ref{L}) we assume that \ $L^{(2)}(\Lambda,z)$ \ is a function of
the similarity variable $V=\frac{\Lambda}{\mu(z)},$ and let $\mathfrak{L}%
(V)=L^{(2)}(\Lambda,z),$ where \ $\mu(z)$ \ is not yet determined. From
(\ref{L}) we get%
\begin{equation}
\frac{\mathfrak{L}^{\prime\prime}}{\mathfrak{L}^{\prime}}=2\Lambda\mu^{\prime
}\frac{\left(  \lambda z+z-\lambda\right)  ^{3}}{\left(  \lambda
z-z-\lambda\right)  \left(  z-\gamma\right)  ^{2}} \label{L1}%
\end{equation}
and (\ref{Linf}) gives%
\begin{equation}
\mathfrak{L}(\infty)=1. \label{Linf1}%
\end{equation}
We can eliminate \ $\Lambda$ by choosing \ $\mu$ \ to satisfy the equation%
\[
2\mu\mu^{\prime}\frac{\left(  \lambda z+z-\lambda\right)  ^{3}}{\left(
\lambda z-z-\lambda\right)  \left(  z-\gamma\right)  ^{2}}=-1
\]
or equivalently%
\begin{equation}
\frac{d\left(  \mu^{2}\right)  }{dz}=-\frac{\left(  \lambda z-z-\lambda
\right)  \left(  z-\gamma\right)  ^{2}}{\left(  \lambda z+z-\lambda\right)
^{3}}. \label{mu1}%
\end{equation}
We choose $\mu(\gamma)=0,$ which is necessary for matching the transition
layer with the corner layer solution (\ref{Fcorner2}), and solve (\ref{mu1})
to obtain%
\begin{equation}
\mu(z)=\sqrt{Y_{2}(z)} \label{mu3}%
\end{equation}
where \ $Y_{2}(z)$ \ was defined in (\ref{Y2}).

Now (\ref{L1}) and (\ref{Linf1}) are just%
\[
\mathfrak{L}^{\prime\prime}=-V\mathfrak{L}^{\prime},\quad\mathfrak{L}%
(\infty)=1
\]
and the solution is%
\[
\mathfrak{L}(V)=\frac{1}{\sqrt{2\pi}}%
{\displaystyle\int\limits_{-\infty}^{V}}
\exp\left(  -\frac{1}{2}\tau^{2}\right)  d\tau.
\]
Thus, the transition layer solution for \ $y-Y_{0}(z)=O\left(  \sqrt
{\varepsilon}\right)  $ and $\gamma<z<1$ is%
\begin{equation}
F_{k}(x)\sim F^{(2)}(V,z)=\sqrt{\varepsilon}\kappa(z)\exp\left[  \frac
{1}{\varepsilon}\Phi(z)\right]  \frac{1}{\sqrt{2\pi}}%
{\displaystyle\int\limits_{-\infty}^{V}}
\exp\left(  -\frac{1}{2}\tau^{2}\right)  d\tau. \label{Ftran}%
\end{equation}

\subsection{Matching the transition layer and $R^{C}$ solutions}

We show that as (\ref{Ftran}) is expanded as $V=\frac{y-Y_{0}(z)}{\mu
(z)\sqrt{\varepsilon}}\rightarrow-\infty,$ the transition layer matches to the
ray expansion, corresponding to rays emanating from $(0,\gamma).$ As
\ $y\rightarrow Y_{0}(z),$ from (\ref{y11}) and (\ref{T0}) we have
\begin{equation}
S(y,z)\sim-\frac{y-Y_{0}(z)}{Y_{2}(z)},\quad T(y,z)\sim\frac{1}{\lambda+1}%
\ln\left(  \frac{z\lambda+z-\lambda}{\rho}\right)  . \label{st}%
\end{equation}
Using (\ref{st}) in (\ref{psi3}) and (\ref{K2}) and expanding for
$y\rightarrow Y_{0}(z),$ we obtain%
\begin{equation}
\Psi\left(  y,z\right)  \sim\Phi(z)-\frac{1}{2}\frac{\left[  y-Y_{0}%
(z)\right]  ^{2}}{Y_{2}(z)} \label{psi4}%
\end{equation}%
\begin{equation}
\mathbb{K}(y,z)\sim-\frac{1}{2\pi}\frac{\sqrt{Y_{2}(z)}}{y-Y_{0}(z)}\frac
{1}{\sqrt{z(1-z)}}. \label{K3}%
\end{equation}

From (\ref{Ftran}) we have
\[
F^{(2)}(V,z)\sim-\sqrt{\varepsilon}\exp\left[  \frac{1}{\varepsilon}%
\Phi(z)\right]  \frac{\kappa(z)}{\sqrt{2\pi}}\frac{1}{V}\exp\left(  -\frac
{1}{2}V^{2}\right)  ,\quad V\rightarrow-\infty
\]
in agreement with (\ref{psi4}) and (\ref{K3}).

\section{The boundary layers at $z=0$ and $z=1$}

\subsection{The boundary layer at $z=0$}

From (\ref{K2}) we see that \ $K(s,t)$ \ is singular as \ $z\rightarrow0.$
Therefore, we find a boundary layer correction near $z=0.$ We consider
solutions of (\ref{diffeq}) which have the asymptotic form%
\begin{equation}
F_{k}(x)-F_{k}(\infty)=F_{k}^{(3)}(y)-F_{k}(\infty)\sim\varepsilon^{\nu_{3}%
-k}\exp\left[  \frac{1}{\varepsilon}\Psi(y,0)\right]  K_{k}^{(3)}(y).
\label{F2}%
\end{equation}
Using (\ref{F2}) in (\ref{diffeq}) and expanding in powers of \ $\varepsilon$
\ gives to leading order%
\begin{equation}
0=\left(  k+1\right)  K_{k+1}^{(3)}+\left[  \gamma\Psi_{y}(y,0)-\lambda
\right]  K_{k}^{(3)}. \label{eq1}%
\end{equation}
Solving (\ref{eq1}) we obtain%
\begin{equation}
K_{k}^{(3)}(y)=\left[  \lambda-\gamma\Psi_{y}(y,0)\right]  ^{k}\frac{1}%
{k!}k^{(3)}(y) \label{K4}%
\end{equation}
and hence
\begin{equation}
F_{k}^{(3)}(y)-F_{k}(\infty)\sim\varepsilon^{\nu_{0}-k}\exp\left[  \frac
{1}{\varepsilon}\Psi(y,0)\right]  \left[  \lambda-\gamma\Psi_{y}(y,0)\right]
^{k}\frac{1}{k!}k^{(3)}(y). \label{F3}%
\end{equation}
Setting \ $k=z/\varepsilon,$ $F_{k}^{(3)}(y)-F_{k}(\infty)=G^{(1)}(y,z)$ \ and
letting \ $k\rightarrow\infty,$ \ we get%
\begin{equation}
G^{(1)}(y,z)\sim\varepsilon^{\nu_{0}}\exp\left\{  \frac{1}{\varepsilon}%
\Psi(y,0)+\frac{1}{\varepsilon}z\ln\left[  \frac{\lambda-\gamma\Psi_{y}%
(y,0)}{z}e\right]  \right\}  \frac{\sqrt{\varepsilon}}{\sqrt{2\pi}}\frac
{1}{\sqrt{z}}k^{(3)}(y). \label{F4}%
\end{equation}

From (\ref{G3}) we have%
\begin{gather}
G(y,z)-G(\infty,z)\sim\frac{\varepsilon}{2\pi}\frac{\sqrt{\rho}}{\sqrt{z}%
}\frac{1}{S(y,0)}\exp\left\{  \frac{1}{\varepsilon}\Psi(y,0)+\frac
{1}{\varepsilon}z\ln\left[  \frac{\lambda-\gamma S(y,0)}{z}e\right]  \right\}
\label{G2}\\
\times\sqrt{\frac{\lambda-\gamma S(y,0)}{\left[  \gamma(1-\gamma
)S(y,0)+\rho\right]  \mathbf{J}_{0}(y)}},\quad z\rightarrow0\nonumber
\end{gather}
with%
\[
S(y,0)=\Psi_{y}(y,0),\quad\mathbf{J}_{0}(y)=\left[  \gamma S(y,0)-\lambda
\right]  \frac{\partial y}{\partial s}\left[  S(y,0),T(y,0)\right]
+\gamma\frac{\partial z}{\partial s}\left[  S(y,0),T(y,0)\right]  .
\]
Matching (\ref{F4}) and (\ref{G2}) we conclude that%
\[
\nu_{3}=\frac{1}{2},\quad k^{(3)}(y)=\frac{\sqrt{\rho}}{\sqrt{2\pi}%
S(y,0)}\sqrt{\frac{\lambda-\gamma S(y,0)}{\left[  \gamma(1-\gamma
)S(y,0)+\rho\right]  \mathbf{J}_{0}(y)}}.
\]
Therefore,%
\begin{align}
F_{k}^{(3)}(y)-F_{k}(\infty)  &  \sim\varepsilon^{\frac{1}{2}-k}\exp\left[
\frac{1}{\varepsilon}\Psi(y,0)\right]  \left[  \lambda-\gamma S(y,0)\right]
^{k}\frac{1}{k!}\label{Fz=0}\\
&  \times\frac{\sqrt{\rho}}{\sqrt{2\pi}S(y,0)}\sqrt{\frac{\lambda-\gamma
S(y,0)}{\left[  \gamma(1-\gamma)S(y,0)+\rho\right]  \mathbf{J}_{0}(y)}%
}.\nonumber
\end{align}
We note that the the right side of (\ref{Fz=0}) is negative since
$S(y,0)<S_{0}<0.$

\subsection{The boundary layer at \ $z=1,\quad0<y<Y_{0}(1)$}

From (\ref{K2}) we see that%
\begin{equation}
\mathbb{K}(y,z)\sim\frac{\sqrt{\rho}}{2\pi S(y,1)}\sqrt{\frac{-\left[
1+\left(  1-\gamma\right)  S(y,1)\right]  }{\left[  \rho+\gamma\left(
1-\gamma\right)  S(y,1)\right]  \mathbf{J}_{1}(y)}}\left(  1-z\right)
^{-\frac{1}{2}},\quad z\rightarrow1 \label{K5}%
\end{equation}
where
\[
\mathbf{J}_{1}(y)=\left[  \left(  1-\gamma\right)  S(y,1)+1\right]
\frac{\partial y}{\partial s}\left[  S(y,1),T(y,1)\right]  -(1-\gamma
)\frac{\partial z}{\partial s}\left[  S(y,1),T(y,1)\right]  ,
\]
\ so that\ $\mathbb{K}(y,z)$\ is singular when \ $z=1.$ Therefore, we
introduce the new variable \ $j=N-k$ \ and consider solutions to
(\ref{diffeq}) of the form%
\begin{equation}
F_{k}(x)=F_{j}^{\left(  4\right)  }(y)\sim\varepsilon^{\nu_{4}-j}\exp\left[
\frac{1}{\varepsilon}\Psi(y,1)\right]  K_{j}^{\left(  4\right)  }(y).
\label{F5}%
\end{equation}
Using (\ref{F5}) in (\ref{diffeq}) gives, as \ $\varepsilon\rightarrow0,$%
\[
\left[  1+(1-\gamma)\Psi_{y}(y,1)\right]  K_{j}^{\left(  4\right)  }%
-\lambda\left(  j+1\right)  K_{j+1}^{\left(  4\right)  }=0
\]
which we can solve to obtain%
\[
K_{j}^{\left(  4\right)  }(y)=\left[  \frac{1+(1-\gamma)\Psi_{y}(y,1)}%
{\lambda}\right]  ^{j}\frac{1}{j!}k^{\left(  4\right)  }(y).
\]
Hence,%
\begin{equation}
F_{j}^{\left(  4\right)  }(y)\sim\varepsilon^{\nu_{4}-j}\exp\left[  \frac
{1}{\varepsilon}\Psi(y,1)\right]  \left[  \frac{1+(1-\gamma)\Psi_{y}%
(y,1)}{\lambda}\right]  ^{j}\frac{1}{j!}k^{\left(  4\right)  }(y). \label{F6}%
\end{equation}

From (\ref{psi3}) we have, as $z\rightarrow1,$%
\begin{equation}
\Psi(y,z)\sim\Psi(y,1)+\ln\left[  \frac{1+(1-\gamma)S\left(  y,1\right)
}{\lambda\left(  1-z\right)  }e\right]  \left(  1-z\right)  ,\quad S\left(
y,1\right)  =\Psi_{y}(y,1). \label{psi}%
\end{equation}
Expanding (\ref{F6}) as $j\rightarrow\infty$ using Stirling's formula and
$j=\frac{1-z}{\varepsilon}$ yields%
\begin{equation}
F_{j}^{\left(  4\right)  }(y)\sim\varepsilon^{\nu_{4}}\frac{\sqrt{\varepsilon
}}{\sqrt{2\pi}}\frac{1}{\sqrt{1-z}}k^{\left(  4\right)  }(y)\exp\left\{
\frac{1}{\varepsilon}\Psi(y,1)+\frac{1}{\varepsilon}\ln\left[  \frac
{1+(1-\gamma)\Psi_{y}(y,1)}{\lambda\left(  1-z\right)  }e\right]  \left(
1-z\right)  \right\}  . \label{F7}%
\end{equation}
Matching (\ref{K5}) and (\ref{psi}) with (\ref{F7}) we conclude that%
\[
\nu_{4}=\frac{1}{2},\quad k^{\left(  4\right)  }(y)=\frac{\sqrt{\rho}}%
{\sqrt{2\pi}S(y,1)}\sqrt{\frac{-\left[  1+\left(  1-\gamma\right)
S(y,1)\right]  }{\left[  \rho+\gamma\left(  1-\gamma\right)  S(y,1)\right]
\mathbf{J}_{1}(y)}}.
\]
Therefore,%
\begin{align}
F_{k}(x)  &  =F_{j}^{\left(  4\right)  }(y)\sim\varepsilon^{\frac{1}{2}-j}%
\exp\left[  \frac{1}{\varepsilon}\Psi(y,1)\right]  \left[  \frac
{1+(1-\gamma)S(y,1)}{\lambda}\right]  ^{j}\frac{1}{j!}\label{F8}\\
&  \times\frac{\sqrt{\rho}}{\sqrt{2\pi}S(y,1)}\sqrt{\frac{-\left[  1+\left(
1-\gamma\right)  S(y,1)\right]  }{\left[  \rho+\gamma\left(  1-\gamma\right)
S(y,1)\right]  \mathbf{J}_{1}(y)}}.\nonumber
\end{align}
In the range $0<y<Y_{0}(1),$ we have $S(y,1)>0.$

\subsection{The corner layer near $\left(  Y_{0}(1),1\right)  $}

When \ \
\[
y\rightarrow Y_{0}(1)=\frac{\rho\left[  \ln\left(  \rho\right)  -1\right]
+1}{\left(  \lambda+1\right)  ^{2}}%
\]
$S\left(  y,1\right)  \rightarrow0,$ and (\ref{F8}) is not defined there. We
consider asymptotic solutions of (\ref{diffeq})\ of the form%
\begin{equation}
F_{k}(x)=F_{j}^{\left(  5\right)  }(y)\sim\left(  \frac{\lambda}{1+\lambda
}\right)  ^{N}\frac{N^{j}}{j!}\lambda^{-j}K_{j}^{\left(  5\right)  }(y)
\label{F9}%
\end{equation}
where $j=N-k$ as before. \ Using (\ref{F9}) in (\ref{diffeq}) gives the
following equation for $K_{j}^{\left(  5\right)  }(y)$%
\[
K_{j}^{\left(  5\right)  }-K_{j+1}^{\left(  5\right)  }=0
\]
with solution%
\[
K_{j}^{\left(  5\right)  }(y)=k^{\left(  5\right)  }(y).
\]
Hence,%
\begin{equation}
F_{j}^{\left(  5\right)  }(y)\sim\left(  \frac{\lambda}{1+\lambda}\right)
^{N}\frac{N^{j}}{j!}\lambda^{-j}k^{\left(  5\right)  }(y). \label{F10}%
\end{equation}
To find $k^{\left(  5\right)  }(y)$ we shall match (\ref{F10}) with the
transition layer solution (\ref{Ftran}).

As $z\rightarrow1,$ (\ref{Ftran}) becomes%
\begin{equation}
F^{(2)}\left[  V\left(  y,z\right)  ,z\right]  \sim\sqrt{\varepsilon}%
\kappa(z)\exp\left[  \frac{1}{\varepsilon}\Phi(z)\right]  \frac{1}{\sqrt{2\pi
}}%
{\displaystyle\int\limits_{-\infty}^{V\left(  y,1\right)  }}
\exp\left(  -\frac{1}{2}\tau^{2}\right)  d\tau\label{F1111}%
\end{equation}
with%
\[
V\left(  y,1\right)  =\frac{y-Y_{0}(1)}{Y_{2}(1)}\frac{1}{\sqrt{\varepsilon}%
}.
\]
From (\ref{F10}) we get as $j=\frac{1-z}{\varepsilon}\rightarrow\infty$%
\begin{equation}
F_{j}^{\left(  5\right)  }(y)\sim\sqrt{\varepsilon}\kappa(z)\exp\left[
\frac{1}{\varepsilon}\Phi(z)\right]  k^{\left(  5\right)  }(y). \label{F121}%
\end{equation}
Matching (\ref{F1111}) and (\ref{F121}) gives%
\[
k^{\left(  5\right)  }(y)=\frac{1}{\sqrt{2\pi}}%
{\displaystyle\int\limits_{-\infty}^{V\left(  y,1\right)  }}
\exp\left(  -\frac{1}{2}\tau^{2}\right)  d\tau
\]
and we conclude that%
\begin{equation}
F_{j}^{\left(  5\right)  }(y)\sim\left(  \frac{\lambda}{1+\lambda}\right)
^{N}\left(  \frac{N}{\lambda}\right)  ^{j}\frac{1}{j!}\frac{1}{\sqrt{2\pi}}%
{\displaystyle\int\limits_{-\infty}^{V\left(  y,1\right)  }}
\exp\left(  -\frac{1}{2}\tau^{2}\right)  d\tau\label{F13}%
\end{equation}
where%
\[
V\left(  y,1\right)  =\frac{y-Y_{0}(1)}{\sqrt{\varepsilon}}\left[
\frac{\left(  1-\gamma\right)  \left(  \rho-4\lambda+1\right)  }{\left(
\lambda+1\right)  ^{3}}-2\frac{2\lambda-\gamma+\left(  \gamma-1\right)
\lambda^{2}}{\left(  \lambda+1\right)  ^{4}}\ln(\rho)\right]  ^{-\frac{1}{2}%
}.\quad
\]

\subsection{The boundary layer at $z=1,\quad Y_{0}(1)<y<\infty$}

To find the solution valid in the region \ $1-z=O(\varepsilon),\quad
Y_{0}(1)<y,$ we immediately conclude that\ $F_{k}(x)$ must be of the form%
\begin{equation}
F_{k}(x)-F_{k}(\infty)=F_{j}^{\left(  6\right)  }(y)-\left(  \frac{\lambda
}{1+\lambda}\right)  ^{N}\binom{N}{j}\lambda^{-j}\sim F_{j}^{\left(  4\right)
}(y), \label{F14}%
\end{equation}
where $F_{j}^{\left(  4\right)  }(y)$ \ is given by (\ref{F8}). This solution
matches to (\ref{G3}) and (\ref{F13}) as $y\rightarrow Y_{0}(1)$ and as
$j\rightarrow\infty,$ respectively.

\section{The boundary $x=0$}

For $x=0$ and $k\leq\left\lfloor c\right\rfloor ,$ the values of $F_{k}(0)$
can be computed from the ray expansion, since $F_{k}(0)-F_{k}(\infty
)\sim\varepsilon\mathbb{K}(0,z)\exp\left[  \frac{1}{\varepsilon}%
\Psi(0,z)\right]  $ is well defined. For $x=0$ and $k\geq\left\lfloor
c\right\rfloor +1,$ we have $F_{k}(0)=0$ by (\ref{FBC}). We now examine how
this boundary condition is satisfied by considering the scale $x=O(1)$
($y=O(\varepsilon))$ and constructing a boundary layer correction to the ray
expansion. Note that this part of the boundary is in the region $R^{C}.$

\subsection{The boundary layer at \ $x=0,\quad\gamma<z<1$}

We shall find the solution satisfying the boundary condition (\ref{BC}). This
boundary condition must be applied on the original $x$-scale. From
(\ref{psi3}) we find that $\Psi(y,z)=\Psi^{(7)}(y,z)+o(y)$, as $y\rightarrow
0,$ where%
\begin{align}
\Psi^{(7)}(y,z)  &  =\left(  z-\gamma\right)  \ln\left[  \frac{ye}%
{(z-\gamma)^{2}}\right]  +(z-1)\ln(1-z)-\ln\left(  \lambda+1\right)
+z\ln(\lambda)\label{psi5}\\
&  -\gamma\ln(\gamma)+\frac{\phi y}{z-\gamma}\ln\left[  \frac{\gamma
y}{\left(  z-\gamma\right)  ^{2}}\right]  +\frac{y}{z-\gamma}\left[
\lambda(1-2\gamma)+(\lambda-1)z\right]  .\nonumber
\end{align}
Hence, we shall consider asymptotic solutions of the form%
\begin{equation}
F_{k}(x)=F^{(7)}(x,z)\sim\varepsilon^{\nu_{7}}\exp\left[  \frac{1}%
{\varepsilon}\Psi^{(7)}(\varepsilon x,z)\right]  K^{(7)}(x,z). \label{F26}%
\end{equation}
Using (\ref{F26}) in (\ref{eqG}) and taking into account that $y=\varepsilon
x$ we get, to leading order,%
\begin{equation}
\frac{\partial K^{(7)}}{\partial z}+\frac{x}{z-\gamma}\frac{\partial K^{(7)}%
}{\partial x}+\left[  \frac{2+z-3\gamma}{2(1-z)(z-\gamma)}\right]  K^{(7)}=0.
\label{K4eq}%
\end{equation}

The most general solution to (\ref{K4eq}) is%
\begin{equation}
K^{(7)}(x,z)=\frac{1}{x}\left(  1-z\right)  ^{\frac{3}{2}}k^{(7)}\left(
\Xi\right)  ,\quad\Xi=\frac{x}{z-\gamma}. \label{K41}%
\end{equation}
Hence,%
\begin{equation}
F^{(7)}(x,z)\sim\varepsilon^{\nu_{7}}\exp\left[  \frac{1}{\varepsilon}%
\Psi^{(7)}(\varepsilon x,z)\right]  \frac{1}{x}\left(  1-z\right)  ^{\frac
{3}{2}}k^{(7)}\left(  \Xi\right)  . \label{F166}%
\end{equation}
To find \ $k^{(7)}\left(  \Xi\right)  $ and \ $\nu_{7}$ \ we will match
(\ref{F166}) with the corner layer solution (\ref{Fcorner2}).

Recalling that $l-\alpha=\frac{z-\gamma}{\varepsilon}$ and using the
asymptotic formula for the Bessel functions \cite{AS},
\[
\text{$J$}_{\nu}(z)\sim\frac{1}{\sqrt{2\pi\nu}}\left(  \frac{ez}{2\nu}\right)
^{\nu}%
\]
we get as \ $\varepsilon\rightarrow0$ and $\vartheta$ fixed%
\begin{equation}
J_{\frac{z-\gamma}{\varepsilon}+\frac{\phi}{\vartheta}}\left(  \frac{\beta
}{\vartheta}\right)  \sim\frac{\sqrt{\varepsilon}}{\sqrt{2\pi}}\frac{1}%
{\sqrt{z-\gamma}}\exp\left\{  \left(  \frac{z-\gamma}{\varepsilon}+\frac{\phi
}{\vartheta}\right)  \ln\left[  \frac{\beta e\varepsilon}{2\vartheta\left(
z-\gamma\right)  }\right]  -\frac{\phi}{\vartheta}\right\}  . \label{J1}%
\end{equation}
Using (\ref{J1}) and writing (\ref{Fcorner2}) in terms of $x=\chi\varepsilon$
and $z=\gamma+(l-\alpha)\varepsilon,$ we have%
\begin{gather}
F_{l}^{(1)}(\chi)\sim\frac{\varepsilon}{2\pi}\sqrt{\frac{\rho}{\phi
\gamma(1-\gamma)(z-\gamma)}}\exp\left\{  \frac{1}{\varepsilon}\left[
\gamma\ln\left(  \frac{\lambda}{\gamma}\right)  -(1-\gamma)\ln(1-\gamma
)-\ln(\lambda+1)\right]  \right\} \nonumber\\
\times\exp\left\{  \frac{z-\gamma}{\varepsilon}\ln\left[  \frac{\lambda
(1-\gamma)e\varepsilon}{\left(  z-\gamma\right)  }\right]  \right\}  \frac
{1}{2\pi i}%
{\displaystyle\int\limits_{Br}}
\frac{1}{\vartheta}\Gamma\left(  \frac{\phi}{\vartheta}+1-\alpha\right)
\exp\left\{  \frac{1}{\varepsilon}\left[  x\vartheta-\left(  z-\gamma\right)
\ln\left(  \vartheta\right)  \right]  \right\} \label{F17}\\
\times\exp\left\{  \frac{\phi}{\vartheta}\ln\left[  \frac{\gamma\varepsilon
}{\phi\left(  z-\gamma\right)  }\right]  -\alpha\ln\left(  \frac{\vartheta
}{\phi}\right)  +\frac{2\lambda\left(  1-\gamma\right)  }{\vartheta}\right\}
d\vartheta.\nonumber
\end{gather}

To evaluate (\ref{F17}) asymptotically as $\varepsilon\rightarrow0$ we shall
use the saddle point method. We find that the integrand has a saddle point at
$\vartheta=\frac{1}{\Xi},$ so that%
\begin{gather}
F_{l}^{(1)}(\chi)\sim\left(  \frac{\varepsilon}{2\pi}\right)  ^{\frac{3}{2}%
}\sqrt{\frac{\rho}{\phi\gamma(1-\gamma)}}\exp\left\{  \frac{1}{\varepsilon
}\left[  z\ln\left(  \lambda\right)  +(z-1)\ln(1-\gamma)-\ln(\lambda
+1)-\gamma\ln(\gamma)\right]  \right\} \label{F19}\\
\times\exp\left\{  \frac{z-\gamma}{\varepsilon}\ln\left[  \frac{e^{2}%
\varepsilon\Xi}{z-\gamma}\right]  +\phi\Xi\ln\left[  \frac{\gamma\varepsilon
}{\phi\left(  z-\gamma\right)  }\right]  +\alpha\ln\left(  \phi\Xi\right)
+2\lambda\left(  1-\gamma\right)  \Xi\right\}  \frac{1}{x}\Xi\Gamma\left(
\phi\Xi+1-\alpha\right) \nonumber
\end{gather}
Taking the limit in (\ref{F166}) as $x\rightarrow0,\ z\rightarrow\gamma$ with
$\Xi$ fixed, we obtain%
\begin{gather}
F^{(7)}(x,z)\sim\varepsilon^{\nu_{7}}\exp\left\{  \phi\Xi\ln\left[
\frac{\gamma\varepsilon\Xi}{z-\gamma}\right]  -\rho\Xi\right\}  \frac{1}%
{x}\left(  1-\gamma\right)  ^{\frac{3}{2}}k^{(7)}\left(  \Xi\right)
\label{F188}\\
\times\exp\left\{  \frac{z-\gamma}{\varepsilon}\ln\left[  \frac{e^{2}%
\varepsilon\Xi}{z-\gamma}\right]  +\frac{1}{\varepsilon}\left[  (z-1)\ln
(1-\gamma)-\ln\left(  \lambda+1\right)  +z\ln(\lambda)-\gamma\ln
(\gamma)\right]  \right\}  .\nonumber
\end{gather}

Matching (\ref{F188}) with (\ref{F19}) we have \
\begin{equation}
k^{(7)}\left(  \Xi\right)  =\left(  \frac{1}{2\pi}\right)  ^{\frac{3}{2}}%
\sqrt{\frac{\rho}{\phi\gamma}}\Gamma\left(  \phi\Xi+1-\alpha\right)  \frac
{\Xi}{(1-\gamma)^{2}}\exp\left[  \left(  \alpha-\phi\right)  \Xi\ln\left(
\phi\Xi\right)  +\phi\Xi\right] \nonumber
\end{equation}
and \ $\nu_{7}=\frac{3}{2}.$ Therefore, for $\gamma<z<1,$
\begin{align}
F^{(7)}(x,z)  &  \sim\left(  \frac{\varepsilon}{2\pi}\right)  ^{\frac{3}{2}%
}x^{\frac{z-\gamma}{\varepsilon}+\alpha}\sqrt{\frac{\rho}{\phi\gamma\left(
1-z\right)  }}\frac{1}{z-\gamma}\left(  \frac{\phi}{z-\gamma}\right)
^{\alpha}\Gamma\left(  \frac{\phi x}{z-\gamma}+1-\alpha\right) \nonumber\\
&  \times\exp\left\{  \frac{1}{\varepsilon}\left(  z-\gamma\right)  \ln\left[
\frac{e\varepsilon}{(z-\gamma)^{2}}\right]  +\frac{1}{\varepsilon}\left[
(z-1)\ln(1-z)-\ln\left(  \lambda+1\right)  +z\ln(\lambda)-\gamma\ln
(\gamma)\right]  \right\} \label{F20}\\
&  \times\exp\left\{  \frac{\phi x}{z-\gamma}\ln\left[  \frac{\gamma
\varepsilon}{\phi\left(  z-\gamma\right)  }\right]  +2\lambda\left(
1-\gamma\right)  \frac{x}{z-\gamma}+(\lambda-1)x\right\} \nonumber
\end{align}
Note that from (\ref{F20}) we have $F_{k}(x)=O\left(  x^{k-\left\lfloor
c\right\rfloor }\right)  ,$ as $x\rightarrow0,$ $k\geq\left\lfloor
c\right\rfloor +1.$

\subsection{Matching the boundary layer at \ $x=0,\quad\gamma<z<1$ and the
$R^{C}$ solution}

Writing \ $x=\frac{y}{\varepsilon},$ and using Stirling's formula we have, as
$\varepsilon\rightarrow0$%
\begin{gather*}
\Gamma\left[  \frac{\phi y}{\left(  z-\gamma\right)  \varepsilon}%
+1-\alpha\right]  \sim\frac{\sqrt{2\pi}}{\sqrt{\varepsilon}}\left[  \frac{\phi
y}{\left(  z-\gamma\right)  \varepsilon}\right]  ^{-\alpha}\sqrt{\frac{\phi
y}{\left(  z-\gamma\right)  }}\\
\times\exp\left\{  \frac{\phi y}{\varepsilon\left(  z-\gamma\right)  }%
\ln\left[  \frac{\phi y}{e\left(  z-\gamma\right)  \varepsilon}\right]
\right\}  .
\end{gather*}
Hence, (\ref{F20}) becomes, for $x=\frac{y}{\varepsilon}\rightarrow\infty$%
\begin{align}
F^{(7)}(x,z)  &  \sim\frac{\varepsilon}{2\pi}\sqrt{\frac{\rho}{\gamma\left(
1-z\right)  }}\frac{1}{\left(  z-\gamma\right)  ^{\frac{3}{2}}}\sqrt
{y}\nonumber\\
&  \times\exp\left\{  \frac{1}{\varepsilon}\left(  z-\gamma\right)  \ln\left[
\frac{ye}{(z-\gamma)^{2}}\right]  +\frac{1}{\varepsilon}\left[  (z-1)\ln
(1-z)-\ln\left(  \lambda+1\right)  +z\ln(\lambda)-\gamma\ln(\gamma)\right]
\right\} \label{F21}\\
&  \times\exp\left\{  \frac{\phi y}{\varepsilon\left(  z-\gamma\right)  }%
\ln\left[  \frac{\gamma y}{\left(  z-\gamma\right)  ^{2}}\right]  +\frac
{y}{\varepsilon}\left[  \lambda-1-\frac{\rho}{z-\gamma}\right]  \right\}
.\nonumber
\end{align}

From (\ref{y(s,t)}) and (\ref{z(s,t)}) we get, as $y\rightarrow0$%
\begin{equation}
s\sim\frac{z-\gamma}{y}+\frac{\phi y}{z-\gamma}\ln\left[  \frac{\gamma
y}{\left(  z-\gamma\right)  ^{2}}\right]  +\frac{1}{z-\gamma}\left[
\lambda(2-3\gamma)+\gamma+(\lambda-1)z\right]  . \label{S}%
\end{equation}
Using (\ref{S}) in (\ref{psi3}) and (\ref{K2}) we find that%
\[
K(s,t)\sim\frac{1}{2\pi}\sqrt{\frac{\rho}{\gamma\left(  1-z\right)  }}\sqrt
{y}\left(  z-\gamma\right)  ^{-\frac{3}{2}},\quad\psi(s,t)\sim\Psi^{(7)}(y,z)
\]
in perfect agreement with (\ref{F21}).

\subsection{The corner layer at $(0,1)$}

For $(y,z)$ close to $(0,1),$ we use the variables $x=\frac{y}{\varepsilon}$
and $j=N-k,$ and look for asymptotic solutions of the form%
\begin{equation}
F_{k}(x)=F_{j}^{(8)}(x)\sim\varepsilon^{\nu_{8}-2j}\exp\left[  \frac
{1}{\varepsilon}\Psi^{(8)}(x;\varepsilon)\right]  K_{j}^{(8)}(x) \label{F22}%
\end{equation}
with%
\begin{equation}
\Psi^{(8)}(x;\varepsilon)=\left(  1-\gamma\right)  \ln\left[  \frac
{x\varepsilon e}{(1-\gamma)^{2}}\right]  -\ln\left(  \lambda+1\right)
+\ln(\lambda)-\gamma\ln(\gamma).\nonumber
\end{equation}
Using (\ref{F22}) in (\ref{diffeq}) gives, to leading order,%
\[
x\lambda(j+1)K_{j+1}^{\left(  8\right)  }=(1-\gamma)^{2}K_{j}^{\left(
8\right)  }%
\]
whose solution is%
\[
K_{j}^{\left(  8\right)  }(x)=\left[  \frac{\left(  1-\gamma\right)  ^{2}%
}{\lambda x}\right]  ^{j}\frac{1}{j!}k^{\left(  8\right)  }\left(  x\right)
.
\]
Hence,%
\begin{equation}
F_{j}^{(8)}(x)\sim\exp\left[  \frac{1}{\varepsilon}\Psi^{(8)}(x;\varepsilon
)\right]  \varepsilon^{\nu_{8}-2j}\left[  \frac{\left(  1-\gamma\right)  ^{2}%
}{\lambda x}\right]  ^{j}\frac{1}{j!}k^{\left(  8\right)  }\left(  x\right)  .
\label{F23}%
\end{equation}
As $j\rightarrow\infty$ (\ref{F23}) gives, by Stirling's formula,%
\begin{equation}
F_{j}^{(8)}(x)\sim F_{j}^{(8)}(x)\sim\frac{1}{\sqrt{2\pi j}}\exp\left[
\frac{1}{\varepsilon}\Psi^{(8)}(x;\varepsilon)+j\right]  \varepsilon^{\nu
_{8}-2j}\left[  \frac{\left(  1-\gamma\right)  ^{2}}{\lambda xj}\right]
^{j}\frac{1}{j!}k^{\left(  8\right)  }\left(  x\right)  . \label{F27}%
\end{equation}
We determine $k^{\left(  8\right)  }\left(  x\right)  $ by matching
(\ref{F27}) to the boundary layer expansion in (\ref{F20}).

Writing \ $z=1-j\varepsilon$ \ and letting $z\rightarrow1,$ we obtain from
(\ref{F20})%
\begin{gather}
F^{(7)}(x,z)\sim\varepsilon\left(  \frac{1}{2\pi}\right)  ^{\frac{3}{2}}%
\Gamma\left(  \frac{\phi x}{1-\gamma}+1-\alpha\right)  \left(  \frac{\phi
x}{1-\gamma}\right)  ^{\alpha}\sqrt{\frac{\rho}{\phi\gamma j}}\frac
{1}{1-\gamma}\varepsilon^{-2j}\left[  \frac{\left(  1-\gamma\right)  ^{2}%
}{\lambda xj}\right]  ^{j}\label{F25}\\
\times\exp\left\{  \frac{1}{\varepsilon}\Psi^{(8)}(x;\varepsilon)+\frac{\phi
x}{1-\gamma}\ln\left[  \frac{\varepsilon\gamma}{\phi\left(  1-\gamma\right)
}\right]  +x(3\lambda-1)\right\}  .\nonumber
\end{gather}
By comparing (\ref{F27}) and (\ref{F25}) we find that%
\[
k^{\left(  8\right)  }\left(  x\right)  =\frac{1}{2\pi}\Gamma\left(
\frac{\phi x}{1-\gamma}+1-\alpha\right)  \left(  \frac{\phi x}{1-\gamma
}\right)  ^{\alpha}\sqrt{\frac{\rho}{\phi\gamma}}\frac{1}{1-\gamma}%
\exp\left\{  \frac{\phi x}{1-\gamma}\ln\left[  \frac{\varepsilon\gamma}%
{\phi\left(  1-\gamma\right)  }\right]  +x(3\lambda-1)\right\}
\]
and $\nu_{8}=1.$Therefore,%
\begin{align*}
F_{j}^{(8)}(x)  &  \sim\frac{\varepsilon}{2\pi}\exp\left\{  \frac
{1}{\varepsilon}\Psi^{(8)}(x;\varepsilon)+\frac{\phi x}{1-\gamma}\ln\left[
\frac{\gamma\varepsilon}{\left(  1-\gamma\right)  \phi}\right]  +(3\lambda
-1)x\right\}  \left[  \frac{\left(  1-\gamma\right)  ^{2}}{\lambda x}\right]
^{j}\frac{1}{j!}\\
&  \times\Gamma\left(  \frac{\phi x}{1-\gamma}+1-\alpha\right)  \left(
\frac{\phi x}{1-\gamma}\right)  ^{\alpha}\sqrt{\frac{\rho}{\phi\gamma}}%
\frac{1}{1-\gamma}\varepsilon^{-2j}.
\end{align*}
We can also show that the above, when expanded for $x\rightarrow\infty,$
matches to the boundary layer expansion in (\ref{F8}), valid for $j=O(1)$ and
$0<y<Y_{0}(1).$

\section{The marginal distribution}

We will now find the equilibrium probability that the buffer content exceeds
$x$%
\begin{equation}
\Pr\left[  X(\infty)>x\right]  =M(x)=1-%
{\displaystyle\sum\limits_{k=0}^{N}}
F_{k}(x) \label{M}%
\end{equation}
for various ranges of $x.$ We will compare our results with those obtained
previously by Morrison \cite{Morrison}.

\subsection{Approximation for $x=O(\varepsilon)=O(1/N)$}

In this region we shall use the spectral representation of the corner layer
solution (\ref{Fcorner}), which applies for $x=\varepsilon\chi=O(\varepsilon
)$. Using the generating function%
\[%
{\displaystyle\sum\limits_{j=-\infty}^{\infty}}
J_{j}(x)z^{j}=\exp\left[  \frac{x}{2}\left(  z-\frac{1}{z}\right)  \right]
\]
we obtain%
\[
\exp\left[  \frac{\rho}{\phi}(j+1-\alpha)\right]  =\left(  \sqrt{u_{0}%
}\right)  ^{-(j+1)}%
{\displaystyle\sum\limits_{l=-\infty}^{\infty}}
J_{l-(j+1)}\left[  -\frac{\beta}{\phi}(j+1-\alpha)\right]  \left(  \sqrt
{u_{0}}\right)  ^{l}.
\]
Therefore,
\begin{align*}
M(x)  &  =M^{(1)}(\chi)\sim%
{\displaystyle\sum\limits_{j\geq0}}
a_{j}\exp\left(  -\frac{\phi}{j+1-\alpha}\chi\right)  \exp\left[  \frac{\rho
}{\phi}(j+1-\alpha)\right]  \left(  \sqrt{u_{0}}\right)  ^{j+1}\\
=  &  \sqrt{\varepsilon}\sqrt{\frac{\rho}{\phi}}\kappa(\gamma)\exp\left[
\frac{1}{\varepsilon}\Phi(\gamma)\right]
{\displaystyle\sum\limits_{j\geq0}}
\frac{\left(  -1\right)  ^{j}}{j!}\frac{1}{\left(  j+1-\alpha\right)  }\\
&  \times\exp\left(  -\frac{\phi}{j+1-\alpha}\chi\right)  \exp\left[
\frac{\rho}{\phi}(j+1-\alpha)+\Upsilon(\theta_{j})\right]  \left(  \sqrt
{u_{0}}\right)  ^{j+1-\alpha}.
\end{align*}
From (\ref{Lambda}) we have%
\[
\exp\left[  \Upsilon(\vartheta_{j})\right]  =\left[  -\left(  j+1-\alpha
\right)  \right]  ^{j+1}\exp\left[  -\frac{2\lambda\left(  1-\gamma\right)
}{\phi}(j+1-\alpha)\right]  \left(  \sqrt{u_{0}}\right)  ^{j+1-\alpha}.
\]
Hence,%
\begin{align}
M^{(1)}(\chi)  &  \sim\sqrt{\varepsilon}\sqrt{\frac{\rho}{\phi}}\kappa
(\gamma)\exp\left[  \frac{1}{\varepsilon}\Phi(\gamma)\right] \label{M1}\\
&  \times%
{\displaystyle\sum\limits_{j\geq0}}
\frac{\left(  j+1-\alpha\right)  ^{j}}{j!}\exp\left(  -\frac{\phi}{j+1-\alpha
}\chi\right)  \exp\left[  \left(  \frac{2\rho}{\phi}-1\right)  (j+1-\alpha
)\right]  \left(  u_{0}\right)  ^{j+1-\alpha}.\nonumber
\end{align}
The formula (\ref{M1}) agrees with Morrison's result (4.14) in \cite{Morrison}%
, taking into account the following notational equivalences

\begin{center}
$%
\begin{tabular}
[c]{|c|c|}\hline
Morrison & Dominici-Knessl\\\hline
$\mu$ & $1-\alpha$\\\hline
$r$ & $\frac{\rho}{\gamma\left(  1-\gamma\right)  }=-S_{0}$\\\hline
$\kappa\left(  \gamma\right)  $ & $-\Phi(\gamma)$\\\hline
$f\left(  \gamma\right)  $ & $-\ln\left(  u_{0}\right)  -2\frac{\rho}{\phi}%
.$\\\hline
\end{tabular}
\ \ \ $
\end{center}

\subsection{Approximation for $x=O(\varepsilon^{-1})=O(N)$}

We shall now use the asymptotic solution in the region $R,$ as given by
(\ref{G3}). We have
\[
M(x)=M^{(2)}(y)=-%
{\displaystyle\sum\limits_{k=0}^{N}}
G\left(  y,\frac{k}{N}\right)  \sim-%
{\displaystyle\int\limits_{0}^{1}}
\exp\left[  \frac{1}{\varepsilon}\Psi(y,z)\right]  \mathbb{K(}y,z)dz
\]
and using the saddle point method we get%
\[
M^{(2)}(y)\sim-\sqrt{\varepsilon}\frac{\sqrt{2\pi}}{\sqrt{-\Psi_{zz}%
(y,\gamma)}}\exp\left[  \frac{1}{\varepsilon}\Psi(y,\gamma)\right]
\mathbb{K(}y,\gamma).
\]
We recall that for a fixed $y$, $\Psi(y,z)$ is maximal at $z=\gamma.$ From
(\ref{U+-}) and (\ref{Tg}) we get%
\[
\Psi_{zz}(y,\gamma)=\frac{1}{\rho}S(y,\gamma),\quad\mathbb{K(}y,\gamma
)=\frac{1}{2\pi S(y,\gamma)}\sqrt{\frac{S(y,\gamma)}{-y_{s}\left(
s,T_{\gamma}\right)  \left[  \gamma\left(  \gamma-1\right)  S(y,\gamma
)-\rho\right]  }}%
\]
where $y_{s}\left(  s,T_{\gamma}\right)  $ is understood to be evaluated at
$s=S(y,\gamma)<S_{0}<0.$ Thus,%
\begin{equation}
M^{(2)}(y)\sim-\frac{\sqrt{\varepsilon}}{\sqrt{2\pi}}\frac{1}{S(y,\gamma
)}\sqrt{\frac{-S_{0}}{y_{s}\left(  s,T_{\gamma}\right)  \left[  S_{0}%
-S(y,\gamma)\right]  }}\exp\left[  \frac{1}{\varepsilon}\Psi(y,\gamma)\right]
. \label{M2}%
\end{equation}
From (\ref{y(s,t)}) and (\ref{Tg}) we get%
\[
y\left(  s,T_{\gamma}\right)  =-\frac{\left[  \phi s+\rho(\lambda+1)\right]
T_{\gamma}+2\rho}{\Delta^{2}}%
\]
and from (\ref{psi3}) and (\ref{Tg}) we have%
\[
\psi(s,T_{\gamma})=sy\left(  s,T_{\gamma}\right)  -\ln(\lambda+1)+\frac{1}%
{2}\left[  \left(  2\gamma-1\right)  s-(\lambda+1)\right]  T_{\gamma}+\frac
{1}{2}\ln\left[  \frac{\lambda s}{\rho+s\gamma(1-\gamma)}\right]  .
\]
The results above agree with Morrison's (5.15) in \cite{Morrison} if we
reconcile notation as below

\begin{center}
$%
\begin{tabular}
[c]{|c|c|}\hline
Morrison & Dominici-Knessl\\\hline
$\tau$ & $-s$\\\hline
$Z(\tau)$ & $y\left(  s,T_{\gamma}\right)  $\\\hline
$\ln\left[  Y(\sigma)\right]  $ & $\Delta(-s)T_{\gamma}(-s)$\\\hline
$U\left(  \tau\right)  $ & $sy\left(  s,T_{\gamma}\right)  -\psi(s,T_{\gamma
}).$\\\hline
\end{tabular}
\ \ \ $
\end{center}

\section{Summary and discussion}

In most of the strip $\mathfrak{D=}\left\{  (y,z):y\geq0,\ 0\leq
z\leq1\right\}  ,$ the asymptotic expansion of $F_{k}(x)=G(y,z)$ is given by%

\begin{equation}
G(y,z)\sim\varepsilon\exp\left[  \frac{1}{\varepsilon}\Psi(y,z)\right]
\mathbb{K}(y,z)\text{ \ in }R^{C} \label{RC1}%
\end{equation}
or%
\begin{equation}
G(\infty,z)-G(y,z)\sim-\varepsilon\exp\left[  \frac{1}{\varepsilon}%
\Psi(y,z)\right]  \mathbb{K}(y,z)\text{ \ in }R. \label{R1}%
\end{equation}
If we consider the continuous part of the density, given by%
\[
f_{k}(x)=F_{k}^{\prime}(x)=\varepsilon\frac{\partial G}{\partial y}(y,z),\quad
x>0,
\]
the transition between $R$ and $R^{C}$ disappears, and we have%
\[
f_{k}(x)\sim\varepsilon\Psi_{y}(y,z)\exp\left[  \frac{1}{\varepsilon}%
\Psi(y,z)\right]  \mathbb{K}(y,z)=\varepsilon\exp\left[  \frac{1}{\varepsilon
}\psi(s,t)\right]  sK(s,t),
\]
everywhere in the interior of$\ \mathfrak{D.}$ Note that $\mathbb{K}(y,z)$
becomes infinite along $y=Y_{0}(z)$ (i.e., $s=0),$ but the product $\Psi
_{y}(y,z)\mathbb{K}(y,z)$ remains finite.

The asymptotic expansion of the boundary probabilities $F_{k}(0),$
$k\leq\left\lfloor c\right\rfloor $ can be obtained by setting $y=0$ in
(\ref{R1}). This expression can be used to estimate the difference
\[
F_{k}(\infty)-F_{k}(0)=\Pr\left[  X(\infty)>0,\quad Z(\infty)=k=\frac
{z}{\varepsilon}\right]  ,\quad z<\gamma
\]
which is exponentially small for $\varepsilon\rightarrow0.$ Also, for a fixed
$z\in\lbrack0,\gamma),$ $f_{k}(x)$ is maximal at $x=0$ (see Figure 4).

\begin{figure}[ptb]
\begin{center}
\rotatebox{270} {\resizebox{10cm}{!}{\includegraphics{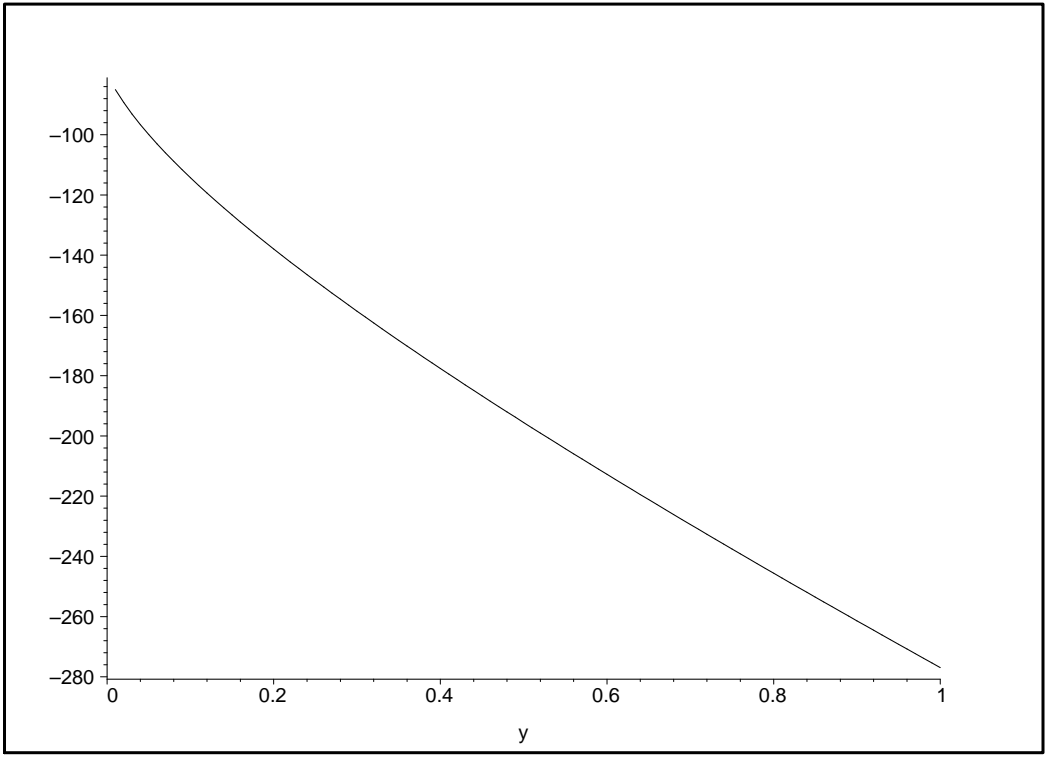}}}
\end{center}
\caption{A plot of $\ln\left\{  \Psi_{y}(y,z)\exp\left[  \frac{1}{\varepsilon
}\Psi(y,z)\right]  \mathbb{K}(y,z)\right\}  $ versus $y$, with $N=100,\ \gamma
=0.2489$ and $z=0.1$.}%
\end{figure}

However, for a fixed $z\in(\gamma,1),$ $f_{k}(x)$ is peeked along the curve
$y=Y_{0}(z)$ (see Figure 5).

\begin{figure}[ptb]
\begin{center}
\rotatebox{270} {\resizebox{10cm}{!}{\includegraphics{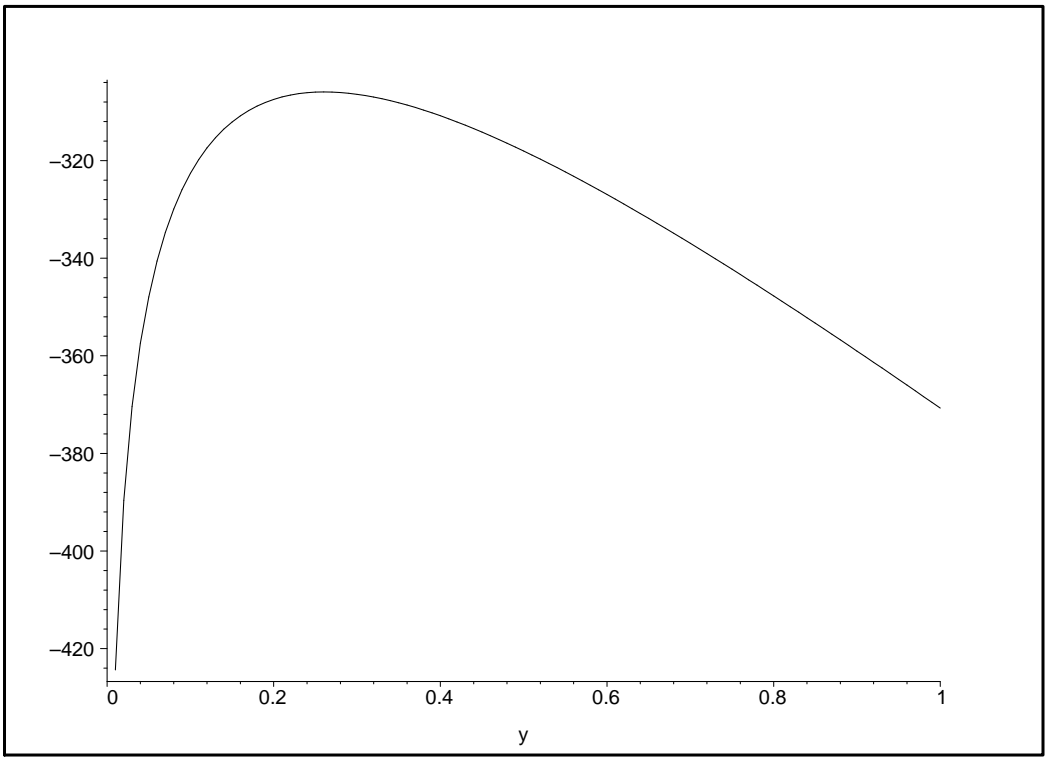}}}
\end{center}
\caption{A plot of $\ln\left\{  \Psi_{y}(y,z)\exp\left[  \frac{1}{\varepsilon
}\Psi(y,z)\right]  \mathbb{K}(y,z)\right\}  $ versus $y$, with $N=100,\ \gamma
=0.2489$ and $z=0.8$.}%
\end{figure}

This means that given $k=zN>c$ active sources, the most likely value of the
buffer will be $x=NY_{0}(z).$ If $zN<c,$ the buffer will most likely be empty.
For a fixed $x\geq0,$ $f_{k}(x)$ achieves its maximum at $z=\gamma$ (see
Figure 6).

\begin{figure}[ptb]
\begin{center}
\rotatebox{270} {\resizebox{10cm}{!}{\includegraphics{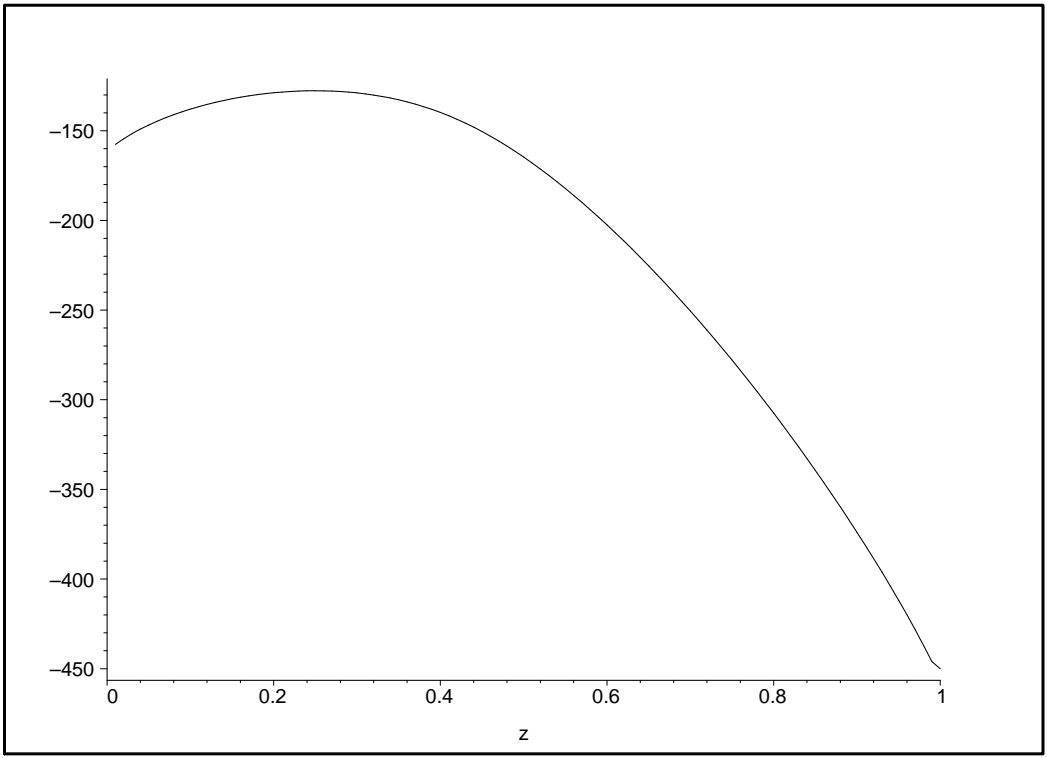}}}
\end{center}
\caption{A plot of $\ln\left\{  \Psi_{y}(y,z)\exp\left[  \frac{1}{\varepsilon
}\Psi(y,z)\right]  \mathbb{K}(y,z)\right\}  $ versus $y$, with $N=100,\ \gamma
=0.2489$ and $y=0.2$.}%
\end{figure}

Below we summarize the various boundary, corner and transition layer
corrections to the results in (\ref{G3}) and (\ref{G4}), where the paragraph
number refers to the corresponding region (see Figure 7).

\begin{figure}[ptbh]
\begin{center}
\rotatebox{0} {\resizebox{15cm}{!}{\includegraphics{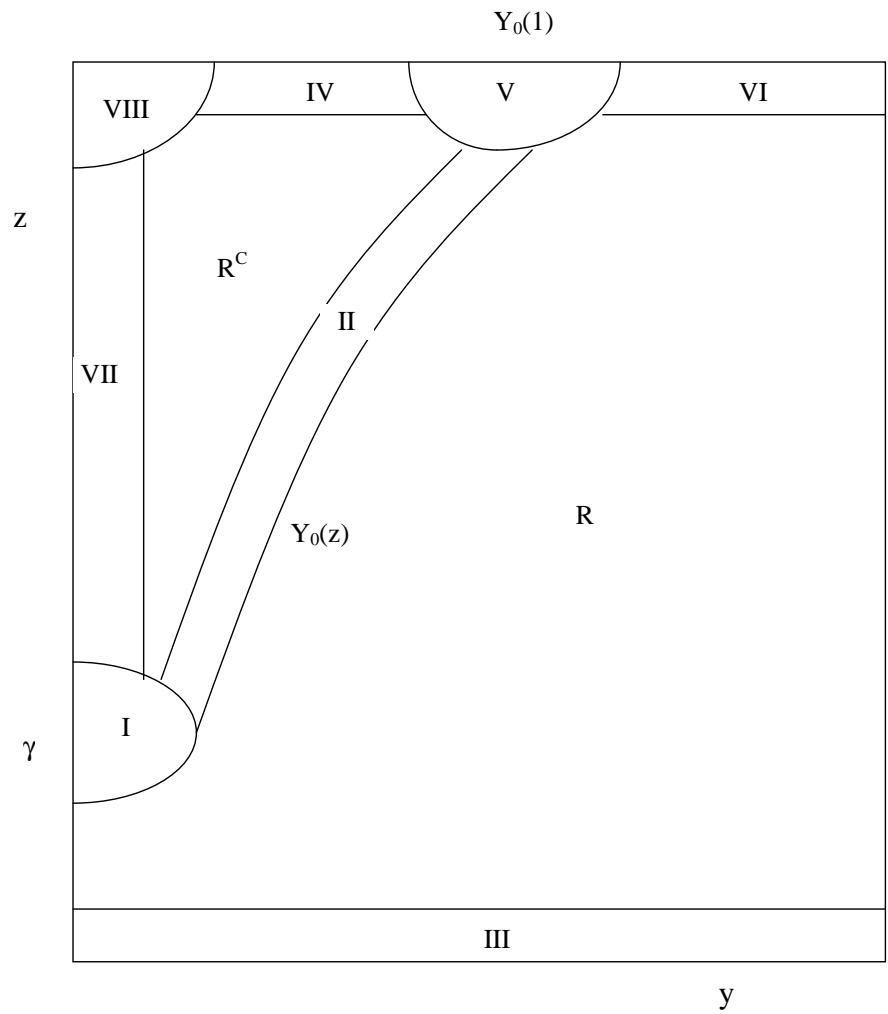}}}
\end{center}
\caption{A sketch of the different asymptotic regions.}%
\end{figure}

\begin{enumerate}
\item $k=l+c-\alpha,$\quad$x=\varepsilon\chi,\quad\chi=O(1)$%
\begin{align*}
F_{k}(x)  &  \sim F_{l}^{(1)}(\chi)=\sqrt{\varepsilon}\sqrt{\frac{\rho}{\phi}%
}\kappa(\gamma)\left(  \sqrt{u_{0}}\right)  ^{l-\alpha}\exp\left[  \frac
{1}{\varepsilon}\Phi(\gamma)\right] \\
&  \times\frac{1}{2\pi i}%
{\displaystyle\int\limits_{\mathrm{Br}}}
e^{\chi\vartheta}\frac{1}{\vartheta}\Gamma\left(  \frac{\phi}{\vartheta
}+1-\alpha\right)  J_{l-\alpha+\frac{\phi}{\vartheta}}\left(  \frac{\beta
}{\vartheta}\right)  \exp\left[  \Upsilon(\vartheta)\right]  d\vartheta.
\end{align*}
where $J_{\cdot}(\cdot)$ denotes the Bessel function, $\mathrm{Br}$ is a
vertical contour in the complex plane with $\operatorname{Re}(\vartheta)>0$
and
\[
\alpha=c-\left\lfloor c\right\rfloor ,\quad\phi=\gamma+\lambda-\gamma
\lambda,\quad\rho=\gamma-\lambda+\lambda\gamma
\]%
\[
\Phi(z)=-z\ln(z)-(1-z)\ln(1-z)+z\ln(\lambda)-\ln(\lambda+1)
\]%
\[
u_{0}=\frac{\lambda}{\gamma}\left(  1-\gamma\right)  ,\quad\kappa(z)=\frac
{1}{\sqrt{2\pi}}\frac{1}{\sqrt{z(1-z)}}%
\]%
\[
\Upsilon(\vartheta)=\left(  \frac{\phi}{\vartheta}-\alpha\right)  \ln\left(
\frac{\vartheta}{\phi}\right)  +\frac{2\lambda\left(  1-\gamma\right)
}{\vartheta}-\frac{\phi}{2\vartheta}\ln\left(  u_{0}\right)  .
\]

\item $y-Y_{0}(z)=O\left(  \sqrt{\varepsilon}\right)  ,\quad\gamma<z<1$%
\[
F_{k}(x)=F^{(2)}(V,z)\sim\sqrt{\varepsilon}\kappa(z)\exp\left[  \frac
{1}{\varepsilon}\Phi(z)\right]  \frac{1}{\sqrt{2\pi}}%
{\displaystyle\int\limits_{-\infty}^{V}}
\exp\left(  -\frac{1}{2}\tau^{2}\right)  d\tau.
\]
with%
\[
V(y,z)=\frac{y-Y_{0}(z)}{\sqrt{\varepsilon}\sqrt{Y_{2}(z)}}%
\]%
\[
Y_{0}(z)=\frac{z-\gamma}{\lambda+1}-\frac{\rho}{\left(  \lambda+1\right)
^{2}}\ln\left(  \frac{z\lambda+z-\lambda}{\rho}\right)  ,\quad\gamma<z<1
\]%
\begin{gather*}
Y_{2}(z)=\frac{2\zeta}{\left(  \lambda+1\right)  ^{4}}\ln\left(
\frac{z+z\lambda-\lambda}{\rho}\right) \\
-\frac{z-\gamma}{\left(  \lambda+1\right)  \left(  \lambda z+z-\lambda\right)
^{2}}\left[  \frac{2\zeta\rho}{\left(  \lambda+1\right)  ^{2}}+\frac{3\zeta
}{\left(  \lambda+1\right)  }\left(  z-\gamma\right)  +\left(  \lambda
-1\right)  \left(  z-\gamma\right)  ^{2}\right]
\end{gather*}%
\[
\zeta=2\lambda-\gamma+(\gamma-1)\lambda^{2}.
\]

\item $k=O(1)$%
\begin{align*}
F_{k}(x)-F_{k}(\infty)  &  =F_{k}^{(3)}(y)-F_{k}(\infty)\sim\varepsilon
^{\frac{1}{2}-k}\exp\left[  \frac{1}{\varepsilon}\Psi(y,0)\right]  \left[
\lambda-\gamma S(y,0)\right]  ^{k}\frac{1}{k!}\\
&  \times\frac{\sqrt{\rho}}{\sqrt{2\pi}S(y,0)}\sqrt{\frac{\lambda-\gamma
S(y,0)}{\left[  \gamma(1-\gamma)S(y,0)+\rho\right]  \mathbf{J}_{0}(y)}}.
\end{align*}

\item $k=N-j,\quad\ j=O(1),\quad0<y<Y_{0}(1)$%
\begin{align*}
F_{k}(x)  &  =F_{j}^{\left(  4\right)  }(y)\sim\varepsilon^{\frac{1}{2}-j}%
\exp\left[  \frac{1}{\varepsilon}\Psi(y,1)\right]  \left[  \frac
{1+(1-\gamma)S(y,1)}{\lambda}\right]  ^{j}\frac{1}{j!}\\
&  \times\frac{\sqrt{\rho}}{\sqrt{2\pi}S(y,1)}\sqrt{\frac{-\left[  1+\left(
1-\gamma\right)  S(y,1)\right]  }{\left[  \rho+\gamma\left(  1-\gamma\right)
S(y,1)\right]  \mathbf{J}_{1}(y)}}.
\end{align*}

\item $k=N-j,\quad\ j=O(1),\quad y-Y_{0}(1)=O\left(  \sqrt{\varepsilon
}\right)  $%
\[
F_{k}(x)=F_{j}^{\left(  5\right)  }(V)\sim\left(  \frac{\lambda}{1+\lambda
}\right)  ^{N}\left(  \frac{N}{\lambda}\right)  ^{j}\frac{1}{j!}\frac{1}%
{\sqrt{2\pi}}%
{\displaystyle\int\limits_{-\infty}^{V}}
\exp\left(  -\frac{1}{2}\tau^{2}\right)  d\tau.
\]

\item $k=N-j,\quad\ j=O(1),\quad y>Y_{0}(1)$%
\[
F_{k}(x)-F_{k}(\infty)=F_{j}^{\left(  6\right)  }(y)-\left(  \frac{\lambda
}{1+\lambda}\right)  ^{N}\binom{N}{j}\lambda^{-j}\sim F_{j}^{\left(  4\right)
}(y),
\]
where $F_{j}^{\left(  4\right)  }(y)$ is as in item 4.

\item $x=O(1),\quad\gamma<z<1$%
\begin{align*}
F_{k}(x)  &  =F^{(7)}(x,z)\sim\left(  \frac{\varepsilon}{2\pi}\right)
^{\frac{3}{2}}\sqrt{\frac{\rho}{\phi\gamma\left(  1-z\right)  }}\frac
{1}{z-\gamma}\left(  \frac{\phi x}{z-\gamma}\right)  ^{\alpha}\Gamma\left(
\frac{\phi x}{z-\gamma}+1-\alpha\right) \\
&  \times\exp\left\{  \frac{1}{\varepsilon}\left(  z-\gamma\right)  \ln\left[
\frac{xe\varepsilon}{(z-\gamma)^{2}}\right]  +\frac{1}{\varepsilon}\left[
(z-1)\ln(1-z)-\ln\left(  \lambda+1\right)  +z\ln(\lambda)-\gamma\ln
(\gamma)\right]  \right\} \\
&  \times\exp\left\{  \frac{\phi x}{z-\gamma}\ln\left[  \frac{\gamma
\varepsilon}{\phi\left(  z-\gamma\right)  }\right]  +2\lambda\left(
1-\gamma\right)  \frac{x}{z-\gamma}+(\lambda-1)x\right\}  .
\end{align*}

\item $k=N-j,\quad\ j=O(1),\quad x=O(1)$%
\begin{align*}
F_{k}(x)  &  =F_{j}^{(8)}(x)\sim\frac{\varepsilon^{1-2j}}{2\pi}\exp\left\{
\frac{1}{\varepsilon}\Psi^{(8)}(x;\varepsilon)+\frac{\phi x}{1-\gamma}%
\ln\left[  \frac{\gamma\varepsilon}{\left(  1-\gamma\right)  \phi}\right]
+(3\lambda-1)x\right\}  \frac{1}{j!}\\
&  \times\left[  \frac{\left(  1-\gamma\right)  ^{2}}{\lambda x}\right]
^{j}\Gamma\left(  \frac{\phi x}{1-\gamma}+1-\alpha\right)  \left(  \frac{\phi
x}{1-\gamma}\right)  ^{\alpha}\sqrt{\frac{\rho}{\phi\gamma}}\frac{1}{1-\gamma
}.
\end{align*}
with%
\[
\Psi^{(8)}(x;\varepsilon)=\left(  1-\gamma\right)  \ln\left[  \frac
{ex\varepsilon}{\left(  1-\gamma\right)  ^{2}}\right]  -\gamma\ln(\gamma
)+\ln\left(  \frac{\lambda}{\lambda+1}\right)  .
\]

\end{enumerate}

\begin{acknowledgement}
The work of C. Knessl was partially supported by NSF Grants DMS 99-71656, DMS
02-02815 and NSA Grant MDA 904-03-1-0036. The work of D. Dominici was
supported in part by NSF Grant 99-73231, provided by Professor Floyd Hanson.
We wish to thank him for his generous sponsorship.
\end{acknowledgement}


\begin{thebibliography}{99}                                                                                               %


\bibitem {AS}M.~Abramowitz and I.~A. Stegun, eds. \emph{Handbook of
mathematical functions with formulas, graphs, and mathematical tables}. Dover
Publications Inc., New York, 1992. Reprint of the 1972 edition.

\bibitem {AMS}D. Anick, D. Mitra, and M. M. Sondhi. \emph{Stochastic theory of
a data-handling system with multiple sources}. Bell System Tech. J.,
\textbf{61} (1982) pp. 1871-1894.

\bibitem {Aziz Liu}A. K. Aziz and J. L. Liu. \emph{A weighted least squares
method for the backward-forward heat equation}. SIAM J. Numer. Anal.,
\textbf{28} (1991), pp. 156-167.

\bibitem {Aziz 2}A. K. Aziz and J. L. Liu. \emph{A Galerkin method for the
forward-backward heat equation}. Math. Comp., \textbf{56} (1991), pp. 35-44.

\bibitem {Baquendi}M. S. Baquendi and P. Grisvard. \emph{Sur une equation d'
evolution changeant de type}. J. Funct. Anal., \textbf{2} (1968), pp. 352--367.

\bibitem {Bardos}C. Bardos, R. E. Caflish and B. Nicolaenko. \emph{The Milne
and Kramers problems for the Boltzmann equation of a hard sphere gas}. Comm.
Pure Appl. Math., \textbf{39} (1986), pp. 323--352.

\bibitem {BenderOrzag}C. M. Bender and S. A. Orzag. \emph{Advanced
mathematical methods for scientists and engineers}. McGraw-Hill. 1978.

\bibitem {Choi Choi 20}B. D. Choi and K. B. Choi. \emph{A Markov modulated
fluid queueing system with strict priority}. Telecomm. Systems, \textbf{9}
(1998), pp. 79--95.

\bibitem {Doorn Sche 22}E. A. van Doorn and W. R. W. Scheinhardt.
\emph{Analysis of birth-death fluid queues}. In B.D. Choi, editor, Proc. KAIST
Applied Mathematics Workshop, pp. 13--29, Taejon, Korea, 1996.

\bibitem {Doorn Sche 23}E. A. van Doorn and W. R. W. Scheinhardt. \emph{A
fluid queue driven by an infinite-state birth-death process}. In V. Ramaswami
and P. E. Wirth, editors, Teletraffic Contributions for the Information Age,
Proc. ITC 15, pp. 465--475, Amsterdam, 1997. Elsevier.

\bibitem {Elwalis Mitra 19}A. Elwalid and D. Mitra. \emph{Analysis,
approximations and admission control of a multi-service multiplexing system
with priorities.} Proc. IEEE INFOCOM '95, 1995, pp. 463--472.

\bibitem {Freidlin weinb}M. Freidlin and H. Weinberger. \emph{On a
backward-forward parabolic equation and its regularization}. J. of Diff. Eq.
\textbf{105} (1993), pp. 264--295.

\bibitem {Gevrey1}M. Gevrey. \emph{Sur certaines equations aux derivees
partielles du type parabolique}. Comptes Rendues (26) \textbf{154} (1912), pp. 1785--1788.

\bibitem {Gevrey 2}M. Gevrey. \emph{Sur les equations aux derivees partielles
du type parabolique}, IV. Journal de Mathematique Pure et Appl. (6)
\textbf{10} (1914), pp. 105--137.

\bibitem {HDL}P. S. Hagan, C. R. Doering and C.D. Levermore. \emph{Mean exit
times for particles driven by weakly colored noise}. SIAM J. Appl. Math.
\textbf{49} (1989), pp. 1480--1513.

\bibitem {Hagan-Okendon}P. S. Hagan and J. R. Ockendon. \emph{Half-range
analysis of a counter-current separator}. J. Math. Anal. and Appl.
\textbf{160} (1991), p. 358--378.

\bibitem {Hopf}E. Hopf. \emph{Mathematical problems of radiative equilibrium}.
Cambridge tracts in mathematics and mathematical physics, no. 31. Cambridge
University Press, 1934.

\bibitem {Keller}J. B. Keller. \emph{Rays, waves and asymptotics}. Bull. Amer.
Math. Soc., \textbf{84} (1978), pp. 727--750.

\bibitem {Keller Wainb}J. B. Keller and H. F. Weinberger. \emph{Boundary and
initial boundary-value problems for separable backward-forward parabolic
problems}. J. Math. Phys., 38B (1997), pp. 4343--4353.

\bibitem {Knessl-Keller}C. Knessl and J. B. Keller. \emph{Ray solution of a
backward-forward parabolic problem for data handling systems}. European J.
Appl. Math., \textbf{11} (2000), pp. 1-12.

\bibitem {Liu Gong 21}Y. Liu and W. Gong. \emph{On fluid queueing system with
strict priority}. Proceedings of IEEE Conference on Decision and Control, 2001.

\bibitem {Knessl-Morrison}C. Knessl and J. A. Morrison. \emph{Heavy traffic
analysis of a data handling system with multiple sources}. SIAM J. Appl.
Math., \textbf{51} (1991), pp. 187--213.

\bibitem {Kobashashi}H. Kobayashi and Q. Ren. \emph{A mathematical theory for
transient analysis of communication networks}. IEICE Trans. Commun., E75-B(12)
(1992), pp. 1266--1276.

\bibitem {Kosten1}L. Kosten. \emph{Stochastic theory of a multi-entry buffer,
part 1}. Delft Progress Report, Series F, \textbf{1} (1974), pp. 10--18.

\bibitem {Kosten2}L. Kosten. \emph{Stochastic theory of a multi-entry buffer,
part 2}. Delft Progress Report, Series F, \textbf{1} (1974), pp. 44--50.

\bibitem {Kultrami}V. G. Kulkarni. \emph{Fluid models for single buffer
systems}. In J.H. Dshalalow, editor, Frontiers in queueing. Models and
Applications in Science and Engineering, pp. 321--338, Boca Raton, Florida,
1997. CRC Press.

\bibitem {Maier}R. S. Maier. \emph{Effective bandwidth of Markov fluids with
occupancy-based admission control}. Proceedings of the 33rd Annual Allerton
Conference on Communication, Control, and Computing. Monticello, Illinois,
Oct. 1995, pp. 766--775.

\bibitem {Mandjes Ridder 24}M. Mandjes and A. Ridder. \emph{A large deviations
analysis of the transient of a queue with many Markov fluid inputs:
approximations and fast simulation}. ACM Transactions on Modeling and Computer
Simulation. \textbf{12} (2002), pp. 1 - 26.

\bibitem {Mitra}D. Mitra. \emph{Stochastic fluid models}. Performance 87, P.
J. Courtois and G. Latouche (editors), Elsevier, (North-Holland), 1988, pp. 39-51.

\bibitem {Mitra 18}D. Mitra. \emph{Stochastic theory of a fluid model of
producers and consumers coupled by a buffer}. Adv. Appl. Prob., 20 (1988), pp. 646--676.

\bibitem {Morrison}J. A. Morrison. \emph{Asymptotic analysis of a
data-handling system with many sources}. SIAM J. Appl. Math., \textbf{49}
(1989), pp. 617-637.

\bibitem {Resing}J. Resing. \emph{Fluid queues and their applications in
telecommunications. }Lecture notes of the S-38.215 Special Course in
Networking Technology, Helsinki University of Technology, Spring 2002.

\bibitem {Scheinhardt}W. Scheinhardt. \emph{Markov-modulated and feedback
fluid models}. Thesis. University of Twente, 1998.

\bibitem {S Weiss}A. Shwartz and A. Weiss. \emph{Large deviations for
performance analysis}. Chapman and Hall, New York, 1995.

\bibitem {Weiss1}A. Weiss. \emph{A new technique for analyzing large traffic
systems}. Advances in Applied Probability, \textbf{18} (1986), pp. 506--532.
\end{thebibliography}
\end{document}